\documentstyle{amsppt}

\loadbold
\loadmsam
\loadmsbm

\def\N{{\Bbb N}}

\def\C{{\Bbb C}}
\def\T{{\Bbb T}}
\def\D{{\Bbb D}}
\def\P{{\Bbb P}}

\def\cL{{\Cal L}}
\def\cH{{\Cal H}}
\def\cF{{\Cal F}}
\def\cG{{\Cal G}}

\topmatter

\title
Five-diagonal matrices and zeros of orthogonal polynomials on the unit circle
\endtitle
\leftheadtext{Five-diagonal matrices and zeros of O.P. on the unit circle}
\rightheadtext{M.J. Cantero, L. Moral, L. Vel\'azquez}

\author
M.J. Cantero, L. Moral, L. Vel\'azquez
\endauthor

\affil
Departamento de Matem\'atica Aplicada,
Universidad de Zaragoza,
Spain
\endaffil

\date
September 2001
\enddate

\thanks
The work of the first and second authors was supported by Direcci\'on General de
Ense\~nanza Superior (DGES) of Spain under grant PB 98-1615.
The work of the last author was supported by CAI, ``Programa Europa de Ayudas a la
Investigaci\'on".
\endthanks





\endtopmatter

\noindent {\bf Abstract}

It is shown that monic orthogonal polynomials on the unit circle
are the cha\-rac\-te\-ris\-tic polynomials of certain
five-diagonal matrices depending on the Schur parameters. This
result is achieved through the study of orthogonal Laurent
polynomials on the unit circle. More precisely, it is a
consequence of the five term recurrence relation obtained for
these orthogonal Laurent polynomials, and the one to one
correspondence established between them and the orthogonal
polynomials on the unit circle. As an application, some results
relating the behavior of the zeros of orthogonal polynomials and
the location of Schur parameters are obtained.

\bigskip


\bigskip

\noindent {\it (1991) AMS Mathematics Subject Classification:} 42C05

\bigskip

\noindent {\it Corresponding author:} Leandro Moral

\hskip 90pt Dpto. Matem\'atica Aplicada

\hskip 90pt Universidad de Zaragoza

\hskip 90pt Pza. San Francisco s/n

\hskip 90pt 50009 Zaragoza (Spain)

\bigskip

\noindent {\it Electronical adress:} lmoral{\@}posta.unizar.es

\bigskip

\noindent {\it Fax:} + 34 976 76 11 25

\vfill\eject

\document

\head
1. Introduction
\endhead

It is very well known that orthogonal polynomials on the real line are the
cha\-rac\-te\-ris\-tic polynomials of the principal submatrices of certain tri-diagonal
infinite matrix called Jacobi matrix [1, 10, 19, 53]. The Jacobi matrix is just the
representation of the multiplication operator in the linear space of real polynomials
when orthogonal polynomials are chosen as a basis. The fact that such polynomials
sa\-tis\-fy a three term recurrence relation is the reason for the tri-diagonality
property of the Jacobi matrix. This property makes the Jacobi matrix simple enough to
obtain properties of the orthogonal polynomials using operator theory techniques [2,
10, 14, 20, 22, 29, 42, 43, 46, 53]. For example, one can get results about their
zeros throughout the spectral analysis of the principal submatrices of the Jacobi
matrix [6, 7, 11, 27, 28, 30, 31, 32, 33, 45, 52].

In the last years, there has been an increasing interest in the
zeros of or\-tho\-go\-nal polynomials on the unit circle due to
their applications in discrete systems analysis, in particular, in
digital signal processing [36, 38, 48, 49]. Unfortunately, only
few things are known about these zeros [3, 4, 5, 15, 17, 44, 47,
50, 51], because the situation in the unit circle is rather more
complicated than in the real line. Orthogonal polynomials on the
unit circle satisfy too a three term recurrence relation that can
be related to a tri-diagonal operator [5, 12, 13], but it does not
provide a spectral representation for their zeros. However, it is
possible to reach such a representation just by computing the
matrix corresponding to the multiplication operator in the linear
space of complex polynomials when orthogonal polynomials are
chosen as a basis. The result is an irreducible Hessenberg matrix
[19], much more complicated than the Jacobi matrix on the real
line. So, this does not seem a so easy way to study properties of
orthogonal polynomials on the unit circle.

The aim of this paper is to improve this situation giving a
five-diagonal matrix representation of orthogonal polynomials on
the unit circle, which yields a spectral interpretation for their
zeros. As we will see, this result comes from the matrix
representation of the multiplication operator in the linear space
of Laurent polynomials, when a suitable basis related to the
orthogonal polynomials is chosen. This matrix representation gives
a spectral interpretation for the zeros of orthogonal polynomials
which is much simpler than the one given by the Hessenberg
matrices. So, it provides an easier way to calculate these zeros
and to study their behavior just using standard methods for
eigenvalue problems of banded matrices.

First of all we will fix some notations.
For an arbitrary finite or infinite matrix $M$, $M^T$ is the transpose matrix of $M$,
and $M^*=\overline M^T$. When $M$ is a square matrix, $M_n$ means the principal
submatrix of $M$ of order $n$ and, as usual, if $M$ is finite, $\det M$ is its
determinant.

In what follows $\P := \C[z]$ is the complex vector space of
polynomials in the variable z with complex coefficients. For $n
\geq 0$, $\P_n := \langle 1,z,\dots,z^n \rangle$ is the
corresponding vector subspace of polynomials with degree less or
equal than $n$, while $\P_{-1} := \{0\}$ is the trivial subspace.
As usual, if $p \in \P_n \backslash \P_{n-1}$, $p^*$ is its
reversed polynomial, defined by $p^*(z) := z^n \overline
p(z^{-1})$. $\Lambda := \C[z,z^{-1}]$ denotes the complex vector
space of Laurent polynomials (L-polynomials) in the variable $z$.
For $m\leq n$, we define the vector subspace $\Lambda_{m,n} :=
\langle z^m,z^{m+1},\dots,z^n \rangle$. Also, for any L-polynomial
$f$, we will consider its substar conjugate defined by
$f_*(z)=\overline f(z^{-1})$. $\T:=\{z\in\C||z|=1\}$ and
$\D:=\{z\in\C||z|<1\}$ are called, respectively, the unit circle
and the open unit disk on the complex plane.

Any hermitian linear functional $\cL$ on $\Lambda$ ($\cL[z^{-n}]=\overline{\cL[z^n]},
\; n=0,1,2,\dots$) defines a sesquilinear functional $(\cdot,\cdot)_\cL
\colon\Lambda\times\Lambda \to \C$ by
$$
(f,g)_\cL:=\cL[f(z) \overline g(z^{-1})], \quad f,g\in\Lambda,
$$
and we say that $f,g\in\Lambda$ are orthogonal with respect to $\cL$ if $(f,g)_\cL=0$.
The hermitian functional $\cL$ is quasi-definite if there exists a sequence of
orthogonal polynomials with respect to $\cL$, that is, a sequence $(p_n)_{n\geq0}$ in
$\P$ satisfying
\itemitem{(I)}
$p_n \in \P_n \backslash \P_{n-1}$,
\itemitem{(II)}
$(p_n,p_m)_\cL = \ell_n \delta_{n,m}, \quad \ell_n \neq 0$.

\noindent The last condition can be replaced equivalently by
\itemitem{(III)}
$(p_n,z^k)_\cL = 0$ if $0 \leq k \leq n-1$,
\item{}
$\quad\, (p_n,z^n)_\cL \neq 0$.

Positive definite hermitian functionals ($\ell_n>0$ for all $n$)
coincide with those given by $\cL[f]:=\int_\C f d\mu$, where $\mu$
is a positive measure with an infinite support lying on $\T$. Due
to this reason, the sequence $(p_n)_{n\geq0}$ is called a sequence
of orthogonal polynomials on the unit circle [54] even in the
general quasi-definite case. In particular, when
$\ell_n=(p_n,p_n)_\cL=\pm1$ for all $n$, we say that
$(p_n)_{n\geq0}$ is a sequence of orthonormal polynomials on the
unit circle.

Given a quasi-definite hermitian functional ${\Cal L}$ on $\Lambda$,
$(\varphi_n)_{n\geq0}$ denotes the unique sequence of orthonormal polynomials with
positive leading coefficients, whereas $(\phi_n)_{n\geq0}$ is the unique sequence
of monic orthogonal polynomials. They are related by
$\varphi_n = \kappa_n \phi_n$, where $\kappa_n := |(\phi_n,\phi_n)_\cL|^{-1/2}$.
Thus, if sg$(\cdot)$ is the sign function,
$(\varphi_n,\varphi_n)_\cL =$ sg$\left((\phi_n,\phi_n)_\cL\right)$.
In what follows we consider $\cL$ normalized by the condition $\cL[1]=1$, so,
$\varphi_0(z)=\phi_0(z)=1$ and $\kappa_0=1$.

It is well known that the sequence $(\phi_n)_{n\geq0}$ is determined by the so called
Schur parameters
$a_n := \phi_n(0)$ through the forward recurrence relation
$$
\eqalign{
& \phi_0(z) = 1, \cr
& \phi_n(z) = z \phi_{n-1}(z) + a_n \phi_{n-1}^*(z), \quad n \geq 1.
}
\tag 1.1
$$
From (1.1) we find that
$$
{(\phi_n,\phi_n)_\cL \over (\phi_{n-1},\phi_{n-1})_\cL} = 1-|a_n|^2, \quad n \geq 1.
\tag 1.2
$$
and, therefore,
$\varepsilon_n:=(\varphi_n,\varphi_n)_\cL \big/ (\varphi_{n-1},\varphi_{n-1})_\cL=$
sg$\left(1-|a_n|^2\right)$.
Hence, apart from the first Schur parameter, that is always $a_0=1$, in the
quasi-definite (positive definite) case it must be  $|a_n| \neq 1$ ($|a_n| < 1$) for
$n\geq1$.  Moreover, any sequence $(a_n)_{n\in\N}$ with this property yields,
through the recurrence (1.1), a sequence of monic orthogonal polynomials on $\T$,
and the associated normalized functional $\cL$ is unique. The parameters
$\{a_k\}_{k=1}^n \subset \C\backslash\T$ allows to construct the finite set of
polynomials $\{\phi_k\}_{k=0}^n$, that is called the finite segment of orthogonal
polynomials associated to $\{a_k\}_{k=1}^n$. Given a finite segment
$\{\phi_k\}_{k=0}^{n-1}$ of orthogonal polynomials, any polynomial with the form
$\phi_n(z) = z \phi_{n-1}(z) + t \phi_{n-1}^*(z), t \in \C\backslash\T$, is called an
extension of $\{\phi_k\}_{k=0}^{n-1}$. Extensions of a finite segment are the possible
candidates to enlarge it.

Relation (1.2) yields
$$
\kappa_n = \prod_{k=1}^n \rho_k^{-1}, \quad n \geq 1,
\tag 1.3
$$
where $\rho_n:=|1-|a_n|^2|^{1/2}$ for $n\geq1$. With this notation we have the
following forward and backward recurrence relations for the orthonormal polynomials
and the reversed ones,
$$
\align
z\varphi_{n-1}(z) & = \rho_n \varphi_n(z) - a_n \varphi_{n-1}^*(z),
\quad n \geq 1,
\tag 1.4 \\
\varphi_{n-1}^*(z) & =  \rho_n \varphi_n^*(z) - \overline a_n z\varphi_{n-1}(z),
\quad n \geq 1,
\tag 1.5 \\
\varphi_n(z) & = a_n \varphi_n^*(z) + \hat\rho_n z\varphi_{n-1}(z),
\quad n \geq 1,
\tag 1.6 \\
\varphi_n^*(z) & = \overline a_n \varphi_n(z) + \hat\rho_n \varphi_{n-1}^*(z),
\quad n \geq 1,
\tag 1.7
\endalign
$$
being $\hat\rho_n := \varepsilon_n \rho_n$ for $n\geq1$.

The finite segment of orthogonal polynomials associated to $\{a_k\}_{k=1}^n$
let us define the $n$-th kernel
$$
K_n(z,y) := \sum_{k=0}^n e_k \varphi_k(z) \overline{\varphi_k(y)},
\tag 1.8
$$
where $e_n:=(\varphi_n,\varphi_n)_\cL$, that is, $e_0=1$ and
$e_n=\prod_{k=1}^n \varepsilon_k$ for $n\geq1$. Using the recurrence relation, the
$n-1$-th kernel can be written equivalently as
$$
K_{n-1}(z,y) =
\cases
\displaystyle{e_n{\varphi_n(z) \overline{\varphi_n(y)}
-\varphi_n^*(z) \overline{\varphi_n^*(y)}
\over z \overline y - 1}}, & $if$ \; z \overline y \neq 1, \\
\displaystyle{e_n z\left(\varphi'_n(z) \overline{\varphi_n(y)}
-(\varphi_n^*)'(z) \overline{\varphi_n^*(y)}\right)}, & $if$ \; z \overline y = 1.
\endcases
\tag 1.9
$$

In what follows it will play an important role the multiplication operator on
$\Lambda$, defined by
$$
\Pi \colon \mathop{\Lambda \longrightarrow \Lambda}\limits_{f(z) \to zf(z)}
$$
Since it lets $\P$ invariant, using the orthonormal polynomials $(\varphi_n)_{n\geq0}$
as a basis of $\P$, it is possible to obtain the matrix representation of the
restriction of $\Pi$ to $\P$, giving the following irreducible Hessenberg matrix
[3, 15, 17, 18, 55]
$$
{\Cal H} =
\pmatrix
d_{0,0} & d_{0,1} & 0       & 0       & 0       & \cdots \\
d_{1,0} & d_{1,1} & d_{1,2} & 0       & 0       & \cdots \\
d_{2,0} & d_{2,1} & d_{2,2} & d_{2,3} & 0       & \cdots \\
\cdots  & \cdots  & \cdots  & \cdots  & \cdots  & \cdots
\endpmatrix
\tag 1.10
$$
where
$$
d_{n,j} =
\cases
-\overline a_j a_{n+1} \prod_{k=j+1}^n \hat\rho_k, & $if$ \; j=0,1,\dots,n-1, \\
-\overline a_n a_{n+1}, & $if$ \; j=n, \\
\hskip 7.5pt \rho_{n+1}, & $if$ \; j=n+1.
\endcases
\tag 1.11
$$
and $(a_n)_{n\geq0}$ are the Schur parameters associated to the orthogonal
polynomials considered.

The characteristic polynomial of the principal submatrix ${\Cal H}_n$ of ${\Cal H}$
of order $n$ is the corresponding $n$-th monic orthogonal polynomial [3, 15, 17, 55].
So, the spectral analysis of ${\Cal H}_n$ can give relations between the zeros of
orthogonal polynomials and the Schur parameters, but, the fact that $\cH$ is a
Hessenberg matrix, together with the complicated dependence of its elements with
respect to the Schur parameters, makes difficult this task. Moreover, the numerical
computations of zeros of high degree orthogonal polynomials, useful, for example, for
applications in digital signal processing, becomes a non trivial problem due to the
Hessenberg structure of $\cH$.

So, a natural question that arises is how to find better spectral
representations for orthogonal polynomials on the unit circle (we
mean with this, the identification of any orthogonal polynomial as
a characteristic polynomial of a matrix). It would be desirable a
situation so similar as possible to the one on the real line, in
particular, we can think on finding banded spectral
representations for orthogonal polynomials on the unit circle.
Moreover, if we want to use this representations to connect the
behavior of zeros of orthogonal polynomials and Schur parameters,
we would need a simple dependence of the elements of the
corresponding matrices with respect to the Schur parameters. We
will find five-diagonal spectral representations for orthogonal
polynomials on the unit circle satisfying this requirement. Apart
from the obvious advantages for the study of orthogonal
polynomials, this result implies a reduction of the eigenvalue
problem for certain Hessenberg matrices to the eigenvalue problem
of a five-diagonal matrix.

The main idea to reach above result is to search inside the matrix
representations of the full operator $\Pi$ on $\Lambda$. This
leads to the study of basis of $\Lambda$ related to orthogonal
polynomials with respect a quasi-definite hermitian functional. As
we will see, the most natural choice are those basis constituted
by Laurent polynomials which are orthogonal with respect to the
same functional.

\head
2. Orthogonal Laurent polynomials on $\T$
\endhead

Orthogonal L-polynomials on the real line appeared in the early eighties in connection
with the theory of continued fractions and strong moment problems [39, 40]. Their
study, not only suffered a rapid development (for a survey, see [37]), but it was
extended to an ampler context, leading to a general theory of rational orthogonal
functions (see [8] and references therein). This theory cover, as a particular case,
the orthogonal polynomials on $\T$. However, the singularities of this particular case
are lost in such a general theory, and, as we will see, these particularities are just
the reason of their utility in the study of orthogonal polynomials on $\T$.

On the real line, orthogonal L-polynomials are a natural generalization of the
orthogonal polynomials when the related functional, initially defined only for
polynomials, is extended to the space of L-polynomials. Our interest is in the
generalization of this idea to the unit circle, where the corresponding functional is
already defined for the full space of L-polynomials.

Although we will deal with orthogonal L-polynomials with respect to a general
quasi-definite hermitian functional on $\T$ (see [23] for the analogous
generalization on the real line), just to understand the following definition, let us
consider a positive definite hermitian functional ${\Cal L}$ on $\Lambda$. Then, the
sesquilinear functional $(\cdot,\cdot)_\cL$ is an inner product on $\Lambda$, and the
orthogonal polynomials with respect to ${\Cal L}$ appear from the standard
orthogonalization of the set $\{1,z,z^2,\dots\}$. Analogously, using the Gram-Schmidt
procedure we can get an orthogonal basis $\{f_n\}_{n\geq0}$ of $\Lambda$ starting from
the ordered basis $\{1,z,z^{-1},z^2,z^{-2},\dots\}$ of $\Lambda$. If we define the
subspaces
$\Lambda_{2n}^+ := \Lambda_{-n,n}, \Lambda_{2n+1}^+ := \Lambda_{-n,n+1}$ for
$n\geq0$ and $\Lambda_{-1}^+ := \{0\}$, then such an orthogonal basis satisfies

\itemitem{(Ia)}
$f_n \in \Lambda_n^+ \backslash \Lambda_{n-1}^+$,
\itemitem{(IIa)}
$(f_n,f_m)_\cL = \ell_n \delta_{n,m}, \quad \ell_n \neq 0$.

So natural than this is to start with the ordered basis
$\{1,z^{-1},z,z^{-2},z^2,\dots\}$. In this situation, the Gram-Schmidt
orthogonalization process gives an orthogonal basis $\{f_n\}_{n\geq0}$ of $\Lambda$
satisfying

\itemitem{(Ib)}
$f_n \in \Lambda_n^- \backslash \Lambda_{n-1}^-$,
\itemitem{(IIb)}
$(f_n,f_m)_\cL = \ell_n \delta_{n,m}, \quad \ell_n \neq 0$,

\noindent where
$\Lambda_{2n}^- := \Lambda_{-n,n}, \Lambda_{2n+1}^- := \Lambda_{-n-1,n}$ for $n\geq0$
and $\Lambda_{-1}^- := \{0\}$.

Above discussion is the origin of the following definition.

\proclaim{Definition 2.1}
A sequence $(f_n)_{n\geq0}$ in $\Lambda$ is called a sequence of right (left)
orthogonal L-polynomials on $\T$ if

\itemitem{(I)}
$f_n \in \Lambda_n^{+(-)} \backslash \Lambda_{n-1}^{+(-)}$

\noindent and there exists a hermitian functional ${\Cal L}$ on $\Lambda$ such that

\itemitem{(II)}
$(f_n,f_m)_\cL = \ell_n \delta_{n,m}, \quad \ell_n \neq 0$.

\noindent Then, we say that $(f_n)_{n\geq0}$ is a sequence of orthogonal
L-polynomials with respect to ${\Cal L}$. If, in addition, $\ell_n = \pm1$ for all
$n$, then we say that $(f_n)_{n\geq0}$ is a sequence of orthonormal L-polynomials.
\endproclaim

\remark{Remark 2.1}
Similarly to orthogonal polynomials, Condition (II) in Definition 2.1 can be replaced
equivalently by
\itemitem{(IIIa)}
$(f_{2n},z^k)_\cL=0$ if $-n+1 \leq k \leq n$,
\quad $(f_{2n},z^{-n})_\cL\neq0$,
\itemitem{}
$(f_{2n+1},z^k)_\cL=0$ if $-n \leq k \leq n$,
\hskip 17pt$(f_{2n+1},z^{n+1})_\cL\neq0$,

\noindent in the case of right orthogonal L-polynomials. For left orthogonal
L-polynomials the equivalent condition is
\itemitem{(IIIb)}
$(f_{2n},z^k)_\cL=0$ if $-n \leq k \leq n-1$,
\quad $(f_{2n},z^n)_\cL\neq0$,
\itemitem{}
$(f_{2n+1},z^k)_\cL=0$ if $-n \leq k \leq n$,
\hskip 17pt$(f_{2n+1},z^{-n-1})_\cL\neq0$.

Analogously to orthogonal polynomials, right and left orthogonal (orthonormal)
L-polynomials are unique up to non null (unimodular) factors.
\endremark

\medskip

Contrary to what happens in the real line, in the unit circle right and left orthogonal
L-polynomials are closely related.

\proclaim{Proposition 2.1}
Let ${\Cal L}$ be a hermitian functional on $\Lambda$ and let $(f_n)_{n\geq0}$ be a
sequence in $\Lambda$. Then, $(f_n)_{n\geq0}$ is a sequence of right orthogonal
(orthonormal) L-polynomials with respect to ${\Cal L}$ iff $(f_{n*})_{n\geq0}$
is a sequence of left orthogonal (orthonormal) polynomials with respect to ${\Cal L}$.
\endproclaim

\demo{Proof}
Obviously $f_{n*} \in \Lambda_n^- \backslash \Lambda_{n-1}^-$ iff
$f_n \in \Lambda_n^+ \backslash \Lambda_{n-1}^+$.
The rest of the proof is just a consequence of the hermiticity of ${\Cal L}$, since it
implies
$(f_{n*},f_{m*})_\cL = \overline{(f_n,f_m)}_\cL$.
$\qed$
\enddemo

Moreover, in the unit circle, orthogonal L-polynomials can be easily constructed from
orthogonal polynomials, something that does not hold in the real line.
This fact, although trivializes such rational functions, is the key for their
usefulness in the study of orthogonal polynomials.

\proclaim{Proposition 2.2}
Let ${\Cal L}$ be a hermitian functional on $\Lambda$ and let $(f_n)_{n\geq0}$ be a
sequence in $\Lambda$. Let us define
$$
\align
& p_{2n}^+(z) = z^n \overline f_{2n}(z^{-1}), \quad p_{2n+1}^+(z) = z^n f_{2n+1}(z),
\quad n \geq 0, \\
& p_{2n}^-(z) = z^n f_{2n}(z), \quad p_{2n+1}^-(z) = z^n \overline f_{2n+1}(z^{-1}),
\quad n \geq 0.
\endalign
$$
Then, $(f_n)_{n\geq0}$ is a sequence of right (left) orthogonal L-polynomials with
respect to ${\Cal L}$ iff $(p_n^{+(-)})_{n\geq0}$ is a sequence of orthogonal
polynomials with respect to ${\Cal L}$. Moreover, $(f_n)_{n\geq0}$ are orthonormal iff
$(p_n^{+(-)})_{n\geq0}$ so are.
\endproclaim

\demo{Proof}
It is straightforward to prove that $p_n^+ \in \P_n \backslash \P_{n-1}$ iff
$f_n \in \Lambda_n^+ \backslash \Lambda_{n-1}^+$. Moreover,
$(p_{2n}^+,z^k)_\cL = \overline{(f_{2n},z^{n-k})}_\cL$ and
$(p_{2n+1}^+,z^k)_\cL = (f_{2n+1},z^{k-n})_\cL$.
Thus, $p_{2n}^+$ is orthogonal to $1,\dots,z^{2n-1}$ iff
$f_{2n}$ is orthogonal to $z^{-n+1},\dots,z^n$, and
$p_{2n+1}^+$ is orthogonal to $1,\dots,z^{2n}$ iff $f_{2n+1}$ is
orthogonal to $z^{-n},\dots,z^n$.
Besides,
$(p_{2n}^+,z^{2n})_\cL = \overline{(f_{2n},z^{-n})}_\cL$ and
$(p_{2n+1}^+,z^{2n+1})_\cL = (f_{2n+1},z^{n+1})_\cL$.
Therefore, according to Remark 2.1,
$(p_n^+)_{n\geq0}$ is a sequence of orthogonal polynomials iff
$(f_n)_{n\geq0}$ is a sequence of right orthogonal L-polynomials.
Finally, since
$(p_n^+,p_n^+)_\cL = (f_n,f_n)_\cL$, we have that
$(p_n^+)_{n\geq0}$ are orthonormal iff $(f_n)_{n\geq0}$ so are.
A similar proof works in the case of left orthogonal L-polynomials.
$\qed$
\enddemo

\remark{Remark 2.2}
Above result establishes in the unit circle a one to one correspondence between
sequences of orthogonal (orthonormal) polynomials and sequences of right or left
orthogonal (orthonormal) L-polynomials. So, any sequence of right orthogonal
(orthonormal) L-polynomials $(f_n)_{n\geq0}$ is obtained from a sequence of orthogonal
(orthonormal) polynomials
$(p_n)_{n\geq0}$ by the relations
$$
\align
& f_{2n}(z) = z^{-n} p_{2n}^*(z), \quad n \geq 0, \\
& f_{2n+1}(z) = z^{-n} p_{2n+1}(z), \quad n \geq 0.
\endalign
$$
The corresponding sequence of left orthogonal (orthonormal) polynomials
$(f_{n*})_{n\geq0}$ is given by
$$
\align
& f_{2n*}(z) = z^{-n} p_{2n}(z), \quad n \geq 0, \\
& f_{2n+1*}(z) = z^{-n-1} p_{2n+1}^*(z), \quad n \geq 0.
\endalign
$$
Moreover,
$\ell_n := (f_n,f_n)_\cL = (f_{n*},f_{n*})_\cL = (p_n,p_n)_\cL$.

Notice that $f_n$ and $f_{n*}$ are linearly independent for $n\geq1$:
for odd index it is obvious;
for even index it is a consequence of the same property for $p_n$ and $p_n^*$.

From above comments we see that the conditions for the existence of orthogonal
polynomials or L-polynomials are the same, that is, quasi-definite hermitian
functionals on $\Lambda$ are just those hermitian functionals for which there exist
(right or left) orthogonal L-polynomials.
\endremark

\medskip

\proclaim{Definition 2.2}
Given a quasi-definite hermitian functional ${\Cal L}$ on $\Lambda$, we denote
by $(\chi_n)_{n\geq0}$ the sequence of right orthonormal L-polynomial
defined by
$$
\align
& \chi_{2n}(z) := z^{-n} \varphi_{2n}^*(z), \quad n \geq 0, \\
& \chi_{2n+1}(z) := z^{-n} \varphi_{2n+1}(z), \quad n \geq 0,
\endalign
$$
where $(\varphi_n)_{n\geq0}$ is the corresponding sequence of orthonormal polynomials
with positive leading coefficients.
We refer to $(\chi_n)_{n\geq0}$ ($(\chi_{n*})_{n\geq0}$) as the standard right
(left) orthogonal L-polynomials associated to ${\Cal L}$.
\endproclaim

The standard right and left orthonormal L-polynomials satisfy some useful relations
that are direct consequences of the recurrence relation for the corresponding
orthonormal polynomials.

\proclaim{Proposition 2.3}
Let ${\Cal L}$ be a quasi-definite hermitian functional on $\Lambda$ and let
$(\chi_n)_{n\geq0}$ be the related sequence of standard right orthonormal
L-polynomials. Then,
$$
\align
& \pmatrix \chi_{2n}(z)  \\ \chi_{2n*}(z)  \endpmatrix
= {1 \over \rho_{2n}} \pmatrix \overline a_{2n} & 1 \\ 1 & a_{2n} \endpmatrix
\pmatrix \chi_{2n-1}(z)  \\ \chi_{2n-1*}(z) \endpmatrix,
\quad n\geq1,
\tag \sl i \\
& \pmatrix \chi_{2n-1}(z) \\ \chi_{2n}(z) \endpmatrix = \Theta_{2n}
\pmatrix \chi_{2n-1*}(z) \\ \chi_{2n*}(z) \endpmatrix,
\quad n\geq1,
\tag \sl ii \\
z & \pmatrix  \chi_{2n*}(z) \\ \chi_{2n+1*}(z) \endpmatrix
= \Theta_{2n+1} \pmatrix \chi_{2n}(z) \\ \chi_{2n+1}(z) \endpmatrix,
\quad n\geq0,
\tag \sl iii \\
& \Theta_n := \pmatrix -a_n & \rho_n \\ \hat\rho_n & \overline a_n \endpmatrix,
\quad n\geq1,
\endalign
$$
where $(a_n)_{n\in\N}$ are the Schur parameters associated to ${\Cal L}$,
$\rho_n = |1-|a_n|^2|^{1/2}$ and $\hat\rho_n = \varepsilon_n \rho_n$ with
$\varepsilon_n =$ {\rm sg}$\left( 1-|a_n|^2 \right)$.
\endproclaim

\demo{Proof} All the relations follow straightforward from (1.4), (1.5), (1.7) and
Definition 2.2.
$\qed$
\enddemo

\subhead
2.1. Recurrence relation for orthogonal L-polynomials on $\T$
\endsubhead

We know that, in the unit circle, the action of the multiplication operator over the
orthogonal polynomials does not provide a recurrence relation for them.
However, as we see in the next proposition, for the orthogonal L-polynomials the
situation is much better.

\proclaim{Proposition 2.4}
Let $(f_n)_{n\geq0}$ be a sequence of right orthogonal L-polynomials on $\T$. Then,
there exist $\pi_{n,k} \in \C$, $|n-k| \leq 2$, such that
$$
\align
& zf_n(z) = \sum_{k=n-2}^{n+2} \pi_{n,k} f_k(z), \quad n \geq 0, \\
& zf_{n*}(z) = \sum_{k=n-2}^{n+2}
{\ell_n \over \ell_k}
\, \pi_{k,n} f_{k*}(z),
\quad n \geq 0,
\endalign
$$
where $\ell_n = (f_n,f_n)_\cL$ and we use the convention $f_k=0$ for $k < 0$.
\endproclaim

\demo{Proof}
Notice that
$zf_n \in z\Lambda_n^+ \subset \Lambda_{n+2}^+ =
\langle f_0,f_1,\dots,f_{n+2} \rangle$.
Moreover, since $f_n$ is orthogonal to
$\langle f_0,f_1,\dots,f_{n-1} \rangle = \Lambda_{n-1}$,
we find that $zf_n$ is orthogonal to $z\Lambda_{n-1} \supset
\Lambda_{n-3} = \langle f_0,f_1,\dots,f_{n-3} \rangle$.
Therefore, $zf_n \in  \langle f_{n-2},\dots,f_{n+2} \rangle$.
Taking into account that $(f_{n*})_{n\geq0}$ is a sequence of left orthogonal
polynomials, we find following similar arguments that
$zf_{n*} \in  \langle f_{n-2*},\dots, f_{n+2*} \rangle$.
Therefore, for $n \geq 0$,
$$
zf_n(z) = \sum_{k=n-2}^{n+2} \pi_{n,k} f_k(z),
\quad zf_{n*}(z) = \sum_{k=n-2}^{n+2} \pi_{n,k*} f_{k*}(z),
$$
with
$\pi_{n,k} = (zf_n,f_k)_\cL \big/ (f_k,f_k)_\cL$ and
$\pi_{n,k*} = (zf_{n*},f_{k*})_\cL \big/ (f_{k*},f_{k*})_\cL$.
Using the definition of the substar conjugate, we find that
$\pi_{n,k*} = \pi_{k,n} \ell_n / \ell_k $.
$\qed$
\enddemo

Above proposition says that orthogonal L-polynomials on the unit circle satisfy a five
term recurrence relation. This fact was already known for orthogonal L-polynomials
on the real line, and used to solve the strong Hamburger moment problem through
operator theory techniques [21].

From the relation between orthogonal L-polynomials and orthogonal
polynomials on $\T$, it is possible to obtain explicitly the
coefficients of the recurrence relation for the orthogonal
L-polynomials in terms of the Schur parameters associated to the
orthogonal polynomials.

\proclaim{Proposition 2.5}
Given a quasi-definite hermitian functional ${\Cal L}$ on $\Lambda$, the related
standard right orthogonal L-polynomials $(\chi_n)_{n\geq0}$ satisfy the recurrence
relation
$$
\align
& z \chi_0(z) = -a_1 \chi_0(z) + \rho_1 \chi_1(z), \\
& z \pmatrix \chi_{2n-1}(z) \\ \chi_{2n}(z) \endpmatrix =
  \widehat M_{2n-1}^T \pmatrix \chi_{2n-2}(z) \\ \chi_{2n-1}(z) \endpmatrix +
  M_{2n} \pmatrix \chi_{2n}(z) \\ \chi_{2n+1}(z) \endpmatrix, \quad n \geq 1, \\
& M_n := \pmatrix -\rho_n a_{n+1} & \rho_n \rho_{n+1} \\
        -\overline a_n a_{n+1} & \overline a_n \rho_{n+1} \endpmatrix,
\quad \widehat M_n := \pmatrix -\hat\rho_n a_{n+1} & \hat\rho_n \hat\rho_{n+1} \\
        -\overline a_n a_{n+1} & \overline a_n \hat\rho_{n+1} \endpmatrix,
\quad n\geq 1,
\endalign
$$
where $(a_n)_{n\in\N}$ are the Schur parameters associated to ${\Cal L}$,
$\rho_n = |1-|a_n|^2|^{1/2}$ and $\hat\rho_n = \varepsilon_n \rho_n$
with $\varepsilon_n =$ {\rm sg}$\left(1-|a_n|^2\right)$.
\endproclaim

\demo{Proof}
From equations (1.4), (1.7) and Definition 2.2, we have that, for $n \geq 1$,
$$
\align
z\chi_{2n-1} & = z^{2-n} \varphi_{2n-1} =
z^{1-n} (\rho_{2n} \varphi_{2n} - a_{2n} \varphi_{2n-1}^*) \\
& = z^{-n} \rho_{2n} (\rho_{2n+1} \varphi_{2n+1} - a_{2n+1} \varphi_{2n}^*)
    - z^{1-n} a_{2n}
    (\overline a_{2n-1} \varphi_{2n-1} + \hat\rho_{2n-1} \varphi_{2n-2}^*) \\
& = \rho_{2n} \rho_{2n+1} \chi_{2n+1} - \rho_{2n} a_{2n+1} \chi_{2n}
    - \overline a_{2n-1} a_{2n} \chi_{2n-1} - \hat\rho_{2n-1} a_{2n} \chi_{2n-2}, \\
z\chi_{2n} & = z^{1-n} \varphi_{2n}^* =
z^{1-n} (\overline a_{2n} \varphi_{2n} + \hat\rho_{2n} \varphi_{2n-1}^*) \\
& = z^{-n} \overline a_{2n} (\rho_{2n+1} \varphi_{2n+1} - a_{2n+1} \varphi_{2n}^*)
    + z^{1-n} \hat\rho_{2n}
    (\overline a_{2n-1} \varphi_{2n-1} + \hat\rho_{2n-1} \varphi_{2n-2}^*) \\
& = \overline a_{2n} \rho_{2n+1} \chi_{2n+1} - \overline a_{2n} a_{2n+1} \chi_{2n}
    + \overline a_{2n-1} \hat\rho_{2n} \chi_{2n-1}
    + \hat\rho_{2n-1} \hat\rho_{2n} \chi_{2n-2}.
\endalign
$$
Besides, from Definition 2.2 and (1.4),
$$
z\chi_0 = z\varphi_0^* = z\varphi_0 = \rho_1 \varphi_1 - a_1 \varphi_0^*
= \rho_1 \chi_1 - a_1 \chi_0,
$$
which completes the proof.
$\qed$
\enddemo

\remark{Remark 2.3}
Following similar arguments to previous proof, from (1.5) and (1.6) we find that the
recurrence relation  for the standard left orthogonal L-polynomials
$(\chi_{n*})_{n\geq0}$ is
$$
\align
& z \pmatrix \chi_{0*}(z) \\ \chi_{1*}(z) \endpmatrix =
  \pmatrix -a_1 \\ \hat\rho_1 \endpmatrix \chi_{0*}(z) +
  M_1 \pmatrix \chi_{1*}(z) \\ \chi_{2*}(z) \endpmatrix, \\
& z \pmatrix \chi_{2n*}(z) \\ \chi_{2n+1*}(z) \endpmatrix =
  \widehat M_{2n}^T
  \pmatrix \chi_{2n-1*}(z) \\ \chi_{2n*}(z) \endpmatrix +
  M_{2n+1} \pmatrix \chi_{2n+1*}(z) \\ \chi_{2n+2*}(z) \endpmatrix,
  \quad n \geq 1.
\endalign
$$
\endremark

\medskip

\head
3. Orthogonal polynomials on $\T$ and five-diagonal matrices
\endhead

The recurrence relation for the standard orthonormal L-polynomials provides a
five-diagonal infinite matrix that plays in the unit circle a similar role to the one
played by the Jacobi matrix in the real line.

\proclaim{Definition 3.1}
The five-diagonal matrix ${\Cal F}$ associated to a quasi-definite hermitian functional
${\Cal L}$ on $\Lambda$ is the following infinite matrix
$$
\align
{\Cal F}
& := \pmatrix
-a_1 \;\, \rho_1 & 0 & 0 & \cdots \\
\quad \widehat M_1^T & M_2 & 0 & \cdots\\
0 & \widehat M_3^T & M_4 & \cdots\\
\vdots & \vdots & \ddots & \ddots
\endpmatrix \\
& = \pmatrix
-a_1 & \rho_1 & 0 \\
-\hat\rho_1 a_2  & \!\!\!\!-\overline a_1 a_2
& \!\!\!-\rho_2 a_3 & \rho_2 \rho_3  \\
\;\;\,\hat\rho_1 \hat\rho_2 & \overline a_1 \hat\rho_2 &
\!\!\!-\overline a_2 a_3 & \overline a_2 \rho_3  & 0 \\
& \!\!\!\!\!\!0 & \!\!\!-\hat\rho_3 a_4  & \!\!\!\!-\overline a_3 a_4
& \!\!\!\!-\rho_4 a_5 & \rho_4 \rho_5 \\
& & \; \hat\rho_3 \hat\rho_4 & \overline a_3 \hat\rho_4
& \!\!\!\!-\overline a_4 a_5 & \overline a_4 \rho_5 & 0 \\
& & & \hskip-30pt\ddots & \hskip-30pt\ddots & \hskip-30pt\ddots & \hskip-15pt\ddots
& \ddots
\endpmatrix
\endalign
$$
where $(a_n)_{n\in\N}$ are the Schur parameters related to ${\Cal L}$,
$\rho_n = |1-|a_n|^2|$ and $\hat\rho_n = \varepsilon_n \rho_n$ with
$\varepsilon_n =$ {\rm sg}$\left(1-|a_n|^2\right)$.
\endproclaim

\remark{Remark 3.1}
Proposition 2.5 means that ${\Cal F}$ is just the matrix of $\Pi$ with respect to
the basis of $\Lambda$ constituted by the standard right orthonormal L-polynomials
$(\chi_n)_{n\geq0}$ related to $\cL$. When $\cL$ is positive definite, Proposition
2.4 implies that $\cF^T$ is the matrix of $\Pi$ when the standard left orthonormal
L-polynomials $(\chi_{n*})_{n\geq0}$ are the basis chosen for $\Lambda$. Taking
into account this proposition we get that, in the general quasi-definite case, the
matrix ${\Cal F}_*$ of $\Pi$ with respect to $(\chi_{n*})_{n\geq0}$ is given by
$\cF_*=E \cF^T E$, being $E$ the infinite diagonal matrix
$$
E:= \pmatrix
e_0 & 0 & 0 & \cdots \\
0 & e_1 & 0 & \cdots \\
0 & 0 & e_2 & \cdots \\
\vdots & \vdots & \vdots & \ddots
\endpmatrix,
$$
with $e_n=(\chi_n,\chi_n)_\cL=(\varphi_n,\varphi_n)_\cL$. In fact,
from Remark 2.3, we find that $\cF_*$ is the result of
substituting in ${\Cal F}^T$ the coefficients $\rho_n$ by
$\hat\rho_n$ and vice versa.
\endremark

Notice that $\Pi = \Pi^{(1)} \Pi^{(2)} = \Pi^{(3)} \Pi^{(1)}$, where $\Pi^{(i)}$,
$i=1,2,3$, are the linear operators on $\Lambda$ defined by
$$
\Pi^{(1)} \colon \mathop{\Lambda \longrightarrow \Lambda}
\limits_{\chi_n(z) \to z\chi_{n*}(z)},
\quad\quad
\Pi^{(2)} \colon \mathop{\Lambda \longrightarrow \Lambda}
\limits_{\chi_{n*}(z) \to \chi_n(z)},
\quad\quad
\Pi^{(3)} \colon \mathop{\Lambda \longrightarrow \Lambda}
\limits_{z\chi_{n*}(z) \to z\chi_n(z)}.
$$
From this fact and Proposition 2.3, we find that
${\Cal F} = {\Cal F}^{(2)} {\Cal F}^{(1)}$ and
${\Cal F}_* = {\Cal F}^{(1)} {\Cal F}^{(2)}$, where ${\Cal F}^{(1)}$,
${\Cal F}^{(2)}$ are the following block-diagonal matrices
$$
{\Cal F}^{(1)} :=
\pmatrix
\Theta_1 & 0 & 0 & \cdots \\
0 & \Theta_3 & 0 & \cdots \\
0 & 0 & \Theta_5 & \cdots \\
\vdots & \vdots & \vdots & \ddots
\endpmatrix,
\quad\quad
{\Cal F}^{(2)} :=
\pmatrix
1 & 0 & 0 & \cdots \\
0 & \Theta_2 & 0 & \cdots \\
0 & 0 & \Theta_4 & \cdots \\
\vdots & \vdots & \vdots & \ddots
\endpmatrix,
$$
being all the blocks two-dimensional, excepting the first one for ${\Cal F}^{(2)}$,
that is one-dimensional.

This gives a decomposition of the five-diagonal matrices ${\Cal F}$, ${\Cal F}_*$ as a
product of two tri-diagonal matrices with special properties:
since $\Theta_n \overline\Theta_n = I_2$ for all $n$, we have that
${\Cal F}^{(i)} \overline{{\Cal F}^{(i)}} = I
= \overline{{\Cal F}^{(i)}} {\Cal F}^{(i)}$, $i=1,2$
($I$ is the infinite unit matrix).
In the positive definite case $\hat\rho_n = \rho_n$ and, thus,
${\Cal F}^{(1)}$, ${\Cal F}^{(2)}$ are symmetric and, so, unitary too.
Therefore, ${\Cal F}$ is unitary for a positive definite functional $\cL$.
This is not a casuality, since, if $\mu$ is the measure on $\T$ associated to the
positive definite functional $\cL$, then ${\Cal F}$ is the matrix of the unitary
operator
$$
U_\mu \colon \mathop{L_\mu^2 \longrightarrow L_\mu^2}
\limits_{f(z) \to zf(z)}
$$
with respect to the Hilbert basis $(\chi_n)_{n\geq0}$ of the space of
$\mu$-square integrable functions $L_\mu^2$.

Notice that, for the principal submatrices of order $n$, we can write too
${\Cal F}_n = {\Cal F}^{(2)}_n {\Cal F}^{(1)}_n$ and
${\Cal F}_{*n} = {\Cal F}^{(1)}_n {\Cal F}^{(2)}_n$,
but now we can only state that
${\Cal F}^{(1)}_n \overline{{\Cal F}^{(1)}_n} = I_n$ for even $n$ and
${\Cal F}^{(2)}_n \overline{{\Cal F}^{(2)}_n} = I_n$ for odd $n$.

\medskip

Analogously to what happens for orthogonal polynomials on $\T$ and
the corresponding Hessenberg matrix, one would expect a relation
between the zeros of orthogonal L-polynomials on $\T$ and the
principal submatrices of the related five-diagonal matrix. But,
the connection between orthogonal polynomials and L-polynomials on
$\T$, provides finally a relation of those principal submatrices
with the zeros of orthogonal polynomials. This fact justifies the
mentioned analogy between the Jacobi matrix on the real line and
the the five-diagonal matrix found on the unit circle.

\proclaim{Theorem 3.1}
Let ${\Cal F}$ be the five-diagonal matrix associated to a quasi-definite hermitian
functional on $\Lambda$ with monic orthogonal polynomials $(\phi_n)_{n\geq0}$ and
standard right orthonormal L-polynomials $(\chi_n)_{n\geq0}$. Then, for $n \geq 1$:

\itemitem{(i)} The characteristic polynomial of the principal submatrix ${\Cal F}_n$ of
${\Cal F}$ of order $n$ is $\phi_n$.

\itemitem{(ii)} The eigenvalues of $\cF_n$ have always geometric multiplicity equal to
1.

\itemitem{(iii)} An eigenvector of $\cF_n$ corresponding to the eigenvalue $\lambda$ is
given by $X_n(\lambda)$, where
$X_n(z):=z^{[{n-1 \over 2}]}(\chi_0(z),\chi_1(z),\dots,\chi_{n-1}(z))^T$.

\endproclaim

\demo{Proof}
From Proposition 2.5, and using 2.3 (iii), we can write
$$
\align
& (zI_n - {\Cal F}_n) X_n(z) = b_n(z), \\
& b_n(z) = \cases
\rho_n z^{1+[{n-1 \over 2}]} \chi_{n*}(z) \; (0,0,\dots,0,1)^T,
& $if$ \; n \; $is even$, \\
\rho_n z^{[{n-1 \over 2}]} \chi_n(z) \; (0,0,\dots,\rho_{n-1},\overline a_{n-1})^T,
& $if$ \; n \; $is odd$,
\endcases
\endalign
$$
where $I_n$ is the unit matrix of order $n$.

If $n$ is even, applying Cramer's rule to solve above system with respect to
$\chi_{n-1}(z)$, we get
$$
\align
\chi_{n-1}(z) = &
{1 \over \det(zI_n-{\Cal F}_n)} \det
\pmatrix
& 0 \\
zI_{n-1}-{\Cal F}_{n-1} & \vdots \\
& 0 \\
\cdots & \rho_n z\chi_{n*}(z) \\
\endpmatrix \\
= & {\det(zI_{n-1}-{\Cal F}_{n-1}) \over \det(zI_n-{\Cal F}_n)}
\rho_n z\chi_{n*}(z).
\endalign
$$
Since $\varphi_j(z)=\kappa_j\phi_j(z)$ and $\rho_j=\kappa_{j-1}/\kappa_j$, using
Definition 2.2 we find that
$$
{\phi_n(z) \over \det(zI_n-{\Cal F}_n)} =
{\phi_{n-1}(z) \over \det(zI_{n-1}-{\Cal F}_{n-1})}.
$$

When $n$ is odd, Cramer's rule to solve the initial system with respect to
$\chi_{n-2}(z)$ gives
$$
\align
\chi_{n-2}(z) = &
{1 \over \det(zI_n-{\Cal F}_n)}
\det \pmatrix
& 0 & 0 \\
zI_{n-2}-{\Cal F}_{n-2} & \vdots & \vdots \\
& 0 & 0 \\
\cdots & \rho_{n-1} \rho_n \chi_n(z) & \rho_{n-1} a_n \\
\cdots & \overline a_{n-1} \rho_n \chi_n(z) & z + \overline a_{n-1} a_n
\endpmatrix \\
= & {\det(zI_{n-2}-{\Cal F}_{n-2}) \over \det(zI_n-{\Cal F}_n)}
z \rho_{n-1} \rho_n \chi_n(z),
\endalign
$$
and, therefore,
$$
{\phi_n(z) \over \det(zI_n-{\Cal F}_n)} =
{\phi_{n-2}(z) \over \det(zI_{n-2}-{\Cal F}_{n-2})}.
$$

Hence, we find by induction that, for $n\geq1$,
$$
{\phi_n(z) \over \det(zI_n-{\Cal F}_n)} =
{\phi_1(z) \over \det(zI_1-{\Cal F}_1)},
$$
which proves (i) since $\phi_1(z) = z + a_1 = \det(zI_1-{\Cal F}_1)$.

So, the eigenvalues of $\cF_n$ coincide with the zeros of $\phi_n$. Using Definition
2.2, we find that
$$
b_n(z) = \cases
\rho_n \varphi_n(z) \; (0,0,\dots,0,1)^T,
& $if$ \; n \; $is even$, \\
\rho_n \varphi_n(z) \; (0,0,\dots,\rho_{n-1},\overline a_{n-1})^T,
& $if$ \; n \; $is odd$.
\endcases
$$
Thus, if $\lambda$ is an eigenvalue of $\cF_n$, it must be
$(\lambda I_n-\cF_n)X_n(\lambda)=0$ and, hence, just showing that
$X_n(\lambda)\neq0$, (iii) is proved. First of all, notice that
$X_n(z)$ is well defined for any value of $z$ since their
components are polynomials in $z$. If $\lambda\neq0$, the first
component of $X_n(\lambda)$ is non null because $\chi_0(z)=1$. On
the contrary, when $\lambda=0$, one of the last two components can
not vanish because
$$
X_n(0)=
\cases
(0,0,\dots,\kappa_{n-2},\kappa_{n-1} a_{n-1})^T, & $if$ \; n \; $is even$, \\
(0,0,\dots,0,\kappa_{n-1})^T, & $if$ \; n \; $is odd$.
\endcases
$$

It only remains to prove (ii).
If $v = (v_1,v_2,\dots,v_n)^T$ satisfies ${\Cal F}_n v = \lambda v$, then
$$
\cases
- a_1 v_1 + \rho_1 v_2 = \lambda v_1, & \\
\widehat M_{2k-1}^T \pmatrix v_{2k-1} \\ v_{2k} \endpmatrix
+ M_{2k} \pmatrix v_{2k+1} \\ v_{2k+2} \endpmatrix
=  \lambda \pmatrix v_{2k} \\ v_{2k+1} \endpmatrix,
& k=1,\dots,\left[{n \over 2}\right]-1, \\
- \hat\rho_{n-1} a_n v_{n-1} - \overline a_{n-1} a_n v_n = \lambda v_n,
& $if$ \; n \; $is even$, \\
\widehat M_{n-2}^T \pmatrix v_{n-2} \\ v_{n-1} \endpmatrix
+ \pmatrix -\rho_{n-1} a_n \\ - \overline a_{n-1} a_n \endpmatrix v_n
=  \lambda \pmatrix v_{n-1} \\ v_n \endpmatrix,
& $if$ \; n \; $is odd$.
\endcases
$$

Let us suppose first that $\lambda\neq0$. In this case we are going to prove by
induction that $v_1=0$ iff $v_n=0$ for all $n$. When $v_1=0$, the first equation of
previous system gives $v_2=0$, while the second equation implies that, if
$v_{2k-1} = v_{2k} = 0$ for some $k \leq [n/2]-1$, then
$$
M_{2k} \pmatrix v_{2k+1} \\ v_{2k+2} \endpmatrix
=  \lambda \pmatrix 0 \\ v_{2k+1} \endpmatrix.
$$
Using the expression of $M_n$ given in Proposition 2.5, we get
from this identity that $v_{2k+1} = v_{2k+2} = 0$. Therefore, if
$v_1=0$, then $v_{2k-1} = v_{2k} = 0$ for $k \leq [n/2]$, what
proves that $v=0$ if $n$ is even.  When $n$ is odd, the last
identity of the system for $v$ gives $v_{2n+1} = 0$ and, so, $v=0$
too.

Hence, if $v = (v_1,v_2,\dots,v_n)^T,v' = (v'_1,v'_2,\dots,v'_n)^T$ are
eigenvectors of ${\Cal F}_n$ with the same eigenvalue $\lambda\neq0$, then
$v_1,v'_1\neq0$. Since $w = v'_1 v - v_1 v'$ satisfies ${\Cal F_n}w = \lambda w$ and
$w_1=0$, it must be $w=0$ and, thus, $v,v'$ are linearly dependent.

Let us consider now the case $\lambda =0$. If $v = (v_1,v_2,\dots,v_n)^T$ satisfies
${\Cal F}_n v = 0$, then, taking into account that $a_n=0$, we have that
$$
\cases
- a_1 v_1 + \rho_1 v_2 = 0, & \\
\widehat M_{2k-1}^T \pmatrix v_{2k-1} \\ v_{2k} \endpmatrix
+ M_{2k} \pmatrix v_{2k+1} \\ v_{2k+2} \endpmatrix = 0,
& k=1,\dots,\left[{n \over 2}\right]-1, \\
\widehat M_{n-2}^T \pmatrix v_{n-2} \\ v_{n-1} \endpmatrix = 0,
& $if$ \; n \; $is odd$.
\endcases
$$
The last equation is equivalent to
$\hat\rho_{n-2} v_{n-2} + \overline a_{n-2} v_{n-1} = 0$,
while, multiplying on the left by $\overline \Theta_{2k}$, the second relation becomes
$\hat\rho_{2k-1} v_{2k-1} + \overline a_{2k-1} v_{2k} =
-a_{2k+1} v_{2k+1} + \rho_{2k+1} v_{2k+2} = 0$.
Therefore, the system for $v$ is equivalent to
$$
\cases
\Theta_{2k-1} \pmatrix v_{2k-1} \\ v_{2k} \endpmatrix = 0,
& k=1,\dots,\left[{n-1\over2}\right], \\
-a_{n-1} v_{n-1} + \rho_{n-1} v_n = 0,
& $if$ \; n \; $is even$.
\endcases
$$
from what it is straightforward to see that $v$ must be proportional to $X_n(0)$.
$\qed$
\enddemo

\remark{Remark 3.2}
From Remark 2.3, it can be shown that (i) and (ii) remains true for the matrix
${\Cal F}_*$. Besides, the eigenvectors associated to an eigenvalue $\lambda$
of $\cF_{*n}$ are given by $X_{n*}(\lambda)$, being $X_{n*}(z):=z^{[{n \over 2}]}
(\chi_{0*}(z),\chi_{1*}(z),\dots,\chi_{n-1*}(z))^T$.
Notice that, when $\lambda=0$,
$$
X_{n*}(0)=
\cases
(0,0,\dots,0,\kappa_{n-1})^T, & $if$ \; n \; $is even$, \\
(0,0,\dots,\kappa_{n-2},\kappa_{n-1} a_{n-1})^T, & $if$ \; n \; $is odd$.
\endcases
$$

Remark 2.1 implies that $\cF_{*n} = E_n \cF^T_n E_n$, where $E_n$ is the
principal submatrix of $E$ of order $n$. Thus, if $\lambda$ is an eigenvalue of
$\cF_n$, and taking into account that $E_n^2=1$, we find that
$\cF^T_n E_n X_{n*}(\lambda) = \lambda E_n X_{n*}(\lambda)$, that
is, $E_n X_{n*}(\lambda)$ is an eigenvector of $\cF^T_n$ with eigenvalue
$\lambda$.

Notice that, if $\lambda$ is a zero of $\phi_n$, then we can take as associated
eigenvectors of $\cF_n$ and $\cF_{*n}$, $V_n(\lambda)$ and $V_{n*}(\lambda)$
respectively, where
$$
\align
& V_n(z)=
\cases
(\chi_0(z),\chi_1(z),\dots,\chi_{n-1}(z))^T,
& $if$ \; \lambda\neq0, \\
(0,0,\dots,\rho_{n-1},a_{n-1})^T, & $if$ \; \lambda=0, \; $even$ \; n,
\\ (0,0,\dots,0,1)^T, & $if$ \; \lambda=0, \; $odd$ \; n,
\endcases \\
& V_{n*}(z)=
\cases
(\chi_{0*}(z),\chi_{1*}(z),\dots,\chi_{n-1*}(z))^T,
& $if$ \; \lambda\neq0, \\
(0,0,\dots,0,1)^T, & $if$ \; \lambda=0, \; $even$ \; n, \\
(0,0,\dots,\rho_{n-1},a_{n-1})^T, & $if$ \; \lambda=0, \; $odd$ \; n.
\endcases
\endalign
$$

Taking into account Remark 3.1, the eigenvalue problem for $\cF_n$ can be
translated into
$$
\cases
{\Cal F}^{(1)}_n X_n(\lambda) = \lambda\overline{{\Cal F}^{(2)}_n} X_n(\lambda),
& $for odd$ \; n, \\
{\Cal F}^{(2)}_n X_{n*}(\lambda) = \lambda\overline{{\Cal F}^{(1)}_n} X_{n*}(\lambda),
& $for even$ \; n.
\endcases
$$

That is, the zeros of $\phi_n$ can be viewed as the eigenvalues of
the five-diagonal matrices $\cF_n$ and $\cF_{*n}$, or,
alternatively, as the generalized eigenvalues of the tri-diagonal
pencil $(\cF^{(1)}_n,\overline{\cF^{(2)}_n})$ or
$(\cF^{(2)}_n,\overline{\cF^{(1)}_n})$ depending on if $n$ is odd
or even.
\endremark

\medskip

Previous theorem gives a spectral interpretation for the zeros of orthogonal
polynomials on $\T$, that allows to calculate them using eigenvalue techniques for
banded matrices. This implies a reduction of their computational cost if
compared with the calculation using Hessenberg matrices [19]. The banded structure
of the five-diagonal matrices, together with their simple dependence on the Schur
parameters, permits even to obtain properties of the zeros by means of standard
matricial techniques, as we will show afterwards. In fact, this banded structure makes
possible to apply similar techniques to those usual for the Jacobi tri-diagonal matrix
on the real line.

Besides, contrary to the Hessenberg matrix ${\Cal H}$ associated to the orthogonal
polynomials on $\T$, every Schur parameter appears in only finitely many elements of
the matrix ${\Cal F}$. This makes easier for $\cF$ than for $\cH$ the analysis of the
effects of perturbations of the sequence of Schur parameters. In particular, every
modification of a finite number of Schur parameters induces a finite dimensional
perturbation of the five-diagonal matrix $\cF$, something that is not true for the
Hessenberg matrix $\cH$. We will take advantage of these facts in the following
section.

As for the relation between $\cF_n$ and $\cH_n$, the geometric
multiplicity of any eigenvalue of an irreducible Hessenberg matrix
is always one [19] and, so, it follows from Theorem 3.1 that
$\cF_n$ and $\cH_n$, not only have the same characteristic
polynomial, but are indeed similar matrices. This fact was not
obvious since $\cF_n$ and $\cH_n$ are matrix representations of
different truncations of the multiplication operator.

\head
4. Applications
\endhead

As a first application of the five-diagonal representation for orthogonal polynomials
on $\T$, we will derive bounds for their zeros in the general quasi-definite case,
where only very few things are known.

\proclaim{Theorem 4.1}
Let ${\Cal L}$ be a quasi-definite hermitian linear functional on $\Lambda$,
$(a_n)_{n\in\N}$ the corresponding sequence of Schur parameters and
$\{z^n_j\}_{j=1}^n$ the zeros of a $n$-th orthogonal polynomial associated with
${\Cal L}$. Then,
$$
R_1 \leq |a_j| \leq R_2, \;\; 1 \leq j \leq n \, \Rightarrow \,
K_1 \leq |z^n_j| \leq K_2, \;\; 1 \leq j \leq n,
$$
where
$K_1=R_1^2+R_2^2-K_2, \; K_2=(R_2+K)^2,
\; K=\max\left\{|1-R_1^2|^{1/2},|1-R_2^2|^{1/2}\right\}.$
\endproclaim

\demo{Proof}
Applying Gershgorin theorem [19, 24] to the matrix $\cF_n$ we find that their
eigenvalues have to lie on a union of disks $D_j, \; j=1,2,\dots,n$, with
centers
$$
c_j = \cases
-a_1 & $if$ \; j=1, \\
-\overline a_{j-1} a_j, & $if$ \; j=2,\dots,n,
\endcases
$$
and radii bounded by
$$
r_j = \cases
\displaystyle \max_{k \leq n} \rho_k, & $if$ \; j=1, \\
\displaystyle (\max_{k \leq n} \rho_k)^2
+ 2 (\max_{k \leq n} \rho_k)(\max_{k \leq n} |a_k|), & $if$ \; j=2,\dots,n.
\endcases
$$
If $R_1 \leq |a_j| \leq R_2$ for $1 \leq j \leq n$, then $R_1 \leq |c_1| \leq R_2$,
$r_1 \leq K$, and $R_1^2 \leq |c_j| \leq R_2^2$, $r_j \leq K^2 + 2KR_2$, for
$2 \leq j \leq n$. Since $|c_j|-r_j \leq |z| \leq |c_j|+r_j$ for $z \in D_j$, any
eigenvalue $\lambda$ of $\cF_n$ must satisfy
$\min\{K_1,R_1-K\} \leq |\lambda| \leq \max\{K_2,R_2+K\}$.
The theorem follows from the fact that $K_1 \leq R_1-K$ and $K_2 \geq R_2+K$:
$K_2 \geq R_2+K$ iff $R_2+K \geq 1$, which is true since
$R_2+K \geq R_2+|1-R_2^2|^{1/2} \geq 1$;
if $R_2 \leq 1$ then $R_1 \leq 1$ too and, thus,
$K_1 = R_1^2+R_2^2-K_2 \leq R_1+R_2-(R_2+K) = R_1-K$;
when $R_2 \geq 1$ we have that $K_1 \leq R_1-K$ iff $(2R_2-1)K \geq R_1^2-R_1-K^2$,
which is true since $R_1^2-R_1-K^2 \leq R_1^2-R_1-|1-R_1^2| \leq 0$.
$\qed$
\enddemo

So, bounds for the complete sequence of Schur parameters give uniform bounds for the
zeros of orthogonal polynomials. Notice that, when applying Gershgorin theorem to the
principal submatrices $\cH_n$ of the Hessenberg matrix ${\Cal H}$, we do not get in
general uniform bounds for the zeros because each new row includes more non vanishing
elements. In fact, using (1.10) and (1.11) we would have found that the centers of the
corresponding Gershgorin disks are the same ones given before, but the bounds for the
radii would be now
$$
r_j = \cases
\displaystyle (\max_{k \leq n} \rho_k) & $if$ \; j=1, \\
\displaystyle (\max_{k \leq n} \rho_k)
+ \sum_{r=1}^{j-1} (\max_{k \leq n} \rho_k)^r (\max_{k \leq n} |a_k|)^2,
& $if$ \; j=2,\dots,n,
\endcases
$$
which does not always give uniform bounds for the zeros in the quasi-definite case.

The bounds that Theorem 4.1 gives for the zeros of orthogonal polynomials locate them
in an annulus of radius $K_2-K_1 = R_2^2-R_1^2+4KR_2+2K^2$. Thus, the best bounds
appear when $K$ is close to 0, that is, when $R_1,R_2$ are close to 1.
Since $K_1,K_2 \to 1$ when $R_1,R_2 \to 1$, we have the following immediate
consequence.

\proclaim{Theorem 4.2}
Let ${\Cal L}$ be a quasi-definite hermitian linear functional on $\Lambda$,
$(a_n)_{n\in\N}$ the corresponding sequence of Schur parameters and
$\{z^n_j\}_{k=1}^n$ the zeros of a $n$-th orthogonal polynomial associated with
${\Cal L}$. Then, $\forall \epsilon>0, \, \exists \delta>0$ such that
$$
\big| |a_j| - 1 \big| < \delta, \;\; 1 \leq j \leq n \, \Rightarrow \,
\big| |z^n_j| - 1 \big| < \epsilon, \;\; 1 \leq j \leq n.
$$
One choice that ensures this fact is
$\delta=\sqrt{1+{\epsilon^2 \over 4(1+\epsilon)}}-1$.
\endproclaim

\demo{Proof}
Following the notations given in Theorem 4.1, now
$R_1=1-\delta, R_2=1+\delta,$ with $\delta>0,$ and, thus,
$K=\sqrt{2\delta+\delta^2}$ and $K_2=1+\epsilon, K_1=1-\epsilon+\delta^2$, where
$\epsilon=2K(1+\delta+K)$. As $K \to 0$ when $\delta \to 0$, the implication in the
theorem turns out to be true. The value given there for $\delta$ is just the only
positive solution for the equation
$2\sqrt{2\delta+\delta^2}(1+\delta+\sqrt{2\delta+\delta^2})=\epsilon$.
$\qed$
\enddemo

Roughly speaking, this last result says that, when the Schur
parameters are close to the unit circle, the zeros of orthogonal
polynomials so are. This fact is easy to derive for positive
definite functionals because, in this case, both, Schur parameters
and zeros, lie on the open unit disk and, hence, the relation
$$
a_n=(-1)^n \prod_{j=1}^n z_j^n
$$
implies that $|a_n| < |z_j^n| < 1$  for  $1 \leq j \leq n$.
Theorem 4.2 is the generalization of this property for the quasi-definite case.

In this first application we have exploited the banded structure of the five-diagonal
matrix representation for orthogonal polynomials. Now we will show the advantages of
this representation for the analysis of perturbations.

Let us consider the five-diagonal matrix $\cF$ of Definition 3.1 as a function of the
sequence $(a_n)_{n\in\N}$. Then, its principal matrix of order $n$ becomes a function
$\cF_n(a_1,a_2,\dots,a_n)$ of $n$ complex variables. Given
$\{a_j\}_{j=1}^n\subset\C\backslash\T$, the analysis of the spectrum of
$\cF_n(a_1,a_2,\dots,a_n)$ is equivalent to the study of the zeros of the last
polynomial in the associated finite segment of orthogonal polynomials
$\{\phi_j\}_{j=0}^n$. We are going to analyze the variation of these zeros under
perturbations of the parameters $\{a_j\}_{j=1}^n$.

Let $\{a_j(t)\}_{k=1}^n\subset\C\backslash\T$ be a set of
parameters depending on a real or complex variable $t$ in a
neighborhood of $t=0$. We consider $a_j=a_j(0)$ as unperturbed
parameters and use for them previous notations. For the perturbed
ones $\{a_j(t)\}_{k=1}^n$ we adopt following notations:
$\{\phi_j^t\}_{j=0}^n$ is the perturbed finite segment of
orthogonal polynomials, the corresponding orthonormal polynomials
are given by $\varphi_j^t(z)=\kappa_j(t)\phi_j^t(z)$ with
$\kappa_{j-1}(t)/\kappa_j(t)=\rho_j(t)=|1-|a_j(t)|^2|^{1/2}$, and
$K_n^t$ is the related $n$-th kernel. The associated standard
right and left orthogonal L-polynomials are denoted by $\chi_j^t$
and $\chi_{j*}^t$ respectively. Also, $E_n(t)$ is the diagonal
matrix given by
$$
E_n(t) = \pmatrix
1 & 0 & \cdots & 0 \\
0 & e_1(t) & \cdots & 0 \\
\vdots & \vdots & \ddots & \vdots \\
0 & 0 & \cdots & e_{n-1}(t)
\endpmatrix,
\quad e_j(t) = \prod_{k=1}^j \varepsilon_k(t),
$$
where $\varepsilon_j(t)=$ sg$(1-|a_j(t)|^2)$.
Besides, we write
$F_n(t):=\cF_n(a_1(t),a_2(t),\dots,a_n(t))$.
Finally, $X_n^t$ and $X_{n*}^t$ are the related vector polynomials given in
Theorem 3.1 and Remark 3.2, which yields the eigenvectors of $F_n(t)$ and
$F_{*n}(t) := E_n(t) F_n(t)^T E_n(t)$, respectively.
The same holds for the vector polynomials $V_n^t$ and $V_{n*}^t$.

Our aim is to study the variation of the zeros of $\phi_n^t$ as a
function of the (real or complex) variable $t$ for suitable
perturbations. For this purpose it would be useful a
``Hellmann-Feynman" type theorem [25, 26, 34, 35] for the
eigenvalues of $F_n(t)$. Although $F_n(t)$ is not self-adjoint
and, even, neither normal too, we can find such a result taking
advantage of the relation between the eigenvectors of $F_n(t)$ and
$F_n(t)^T$ [19] given in Remark 3.2. This is the basic idea under
the proof of the following proposition.

\proclaim{Proposition 4.1}
If $a_j(t)$ and $\overline{a_j(t)}$ are differentiable at $t=0$ for $j=1,2,\dots,n$,
and we can express an eigenvalue of $F_n(t):=\cF_n(a_1(t),a_2(t),\dots,a_n(t))$
as a function $\lambda(t)$ differentiable at $t=0$, then
$$
V_{n*}(\lambda)^T E_n F'_n(0) V_n(\lambda) =
\cases
K_{n-1}(\lambda,\overline{\lambda^{-1}}) \; \lambda'(0),
& $if$ \; \lambda = \lambda(0) \neq 0, \\
e_{n-1} a_{n-1} \lambda'(0),
& $if$ \; \lambda = \lambda(0) = 0.
\endcases
$$
\endproclaim

\demo{Proof} Under the hypothesis, $E_n(t)$, $F_n(t)$, as well as
the eigenvectors $X_n^t(\lambda(t))$ and $X_{n*}^t(\lambda(t))$,
are differentiable at $t=0$ (notice that $E_n(t)$ has to be
constant in a neighborhood of $t=0$ since
$\{a_j\}_{j=1}^n\subset\C\backslash\T$). Thus, taking derivatives
at $t=0$ in the identity
$$
X_{n*}^t(\lambda(t))^T E_n(t) F_n(t) X_n^t(\lambda(t)) =
\lambda(t) X_{n*}^t(\lambda(t))^T E_n(t) X_n^t(\lambda(t)),
$$
and using the fact that $F_n X_n(\lambda) = \lambda X_n(\lambda)$ and
$F_n^T E_n X_{n*}(\lambda) =
\lambda E_n X_{n*}(\lambda)$, we get that
$$
X_{n*}(\lambda)^T E_n F'_n(0) X_n(\lambda) =
\lambda'(0) X_{n*}(\lambda)^T E_n X_n(\lambda).
$$
From this equation and the proportionality between the vectors
$X_n(\lambda)$, $X_{n*}(\lambda)$ and
$V_n(\lambda)$, $V_{n*}(\lambda)$,
we find the desired result just noticing that
$$
V_{n*}(\lambda)^T E_n V_n(\lambda) =
\cases
K_{n-1}(\lambda,\overline{\lambda^{-1}}),
& $if$ \; \lambda \neq 0, \\
e_{n-1} a_{n-1},
& $if$ \; \lambda = 0.
\endcases
$$
$\qed$
\enddemo

\remark{Remark 4.1}
From the expression (1.9) for the kernel, and taking into account that $\lambda$ is
a zero of $\varphi_n$, for $\lambda \neq 0$ we get that
$$
K_{n-1}(\lambda,\overline{\lambda^{-1}}) =
e_n \lambda^{1-n} \varphi'_n(\lambda) \varphi_n^*(\lambda).
\tag 4.1
$$
Notice that $\varphi_n$ and $\varphi_n^*$ can never have a common non vanishing zero
since, otherwise, the recurrence relations (1.4) and (1.5) would imply that
$\varphi_j$ and $\varphi_j^*$ have this common zero too for all $j \leq n$, which is
impossible. Thus, $K_{n-1}(\lambda,\overline{\lambda^{-1}}) = 0$ iff the zero $\lambda$
of $\phi_n$ is multiple.

On the other hand, if $\lambda = 0$, then $\phi_n$ has a zero at the origin and, so,
it must be $a_n=0$. In this situation, from the recurrence relation (1.1) we have that
$\phi_n(z)=z\phi_{n-1}(z)$ and, hence, $a_{n-1}=0$ iff the zero $\lambda=0$ of $\phi_n$
is multiple.
\endremark

\medskip

As a first application of previous proposition, we will discuss the variation of the
zeros of a $n$-th orthogonal polynomial under the perturbation of the last Schur
parameter $a_n$, which means the study of the zeros of the extensions of a finite
segment of orthogonal polynomials. This is equivalent to analyze the eigenvalues of
$F(t)=\cF_n(a_1,a_2,\dots,a_{n-1},t)$ as a function of the complex variable $t$.

\proclaim{Theorem 4.3}
Let $\{\phi_j\}_{j=0}^{n-1}$ be the finite segment of orthogonal polynomials associated
to $\{a_j\}_{j=1}^{n-1}\subset\C\backslash\T$.
Then, the zeros of $\phi_n^t(z)=z\phi_{n-1}(z)+t\phi_{n-1}^*(z)$ are
simple for $t\in\C\backslash S$, where $S$ has only a finite number of points in any
compact subset of $\C$.  For $t$ in any simply connected domain of $\C\backslash S$,
the zeros of $\phi_n^t$ can be expressed as holomorphic functions of $t$. If
$\lambda(t)$ is one of such functions, then,
$$
\lambda'(t) =
\cases
\displaystyle{
-e_{n-1} \lambda(t)^{1-n}
{(\varphi_{n-1}^*(\lambda(t)))^2 \over
K_{n-1}(\lambda(t),\overline{\lambda(t)^{-1}})},} & $if$ \; \lambda(t) \neq 0, \\
\displaystyle{
-e_{n-1} {1 \over a_{n-1}},}
& $if$ \; \lambda(t) = 0.
\endcases
$$
\endproclaim

\demo{Proof}
Since $F(t):=\cF_n(a_1,a_2,\dots.a_{n-1},t)$ is an analytic function of $t$ in $\C$,
the number of distinct eigenvalues of $F_n(t)$ and their algebraic multiplicities are
constant up to, at most, a set $S$ with only a finite number of points in any compact
subset of $\C$. Moreover, for $t$ in any simply connected domain of $\C\backslash S$,
these eigenvalues can be expressed as holomorphic functions $\lambda(t)$ [41].

So, just applying Proposition 4.1, and simplifying the result
using Proposition 2.3 and the relation between standard orthogonal
L-polynomials and orthogonal polynomials, we get that
$$
\cases
K_{n-1}(\lambda(t),\overline{\lambda(t)^{-1}}) \lambda'(t) =
-e_{n-1} \lambda(t)^{1-n} (\varphi_{n-1}^*(\lambda(t)))^2,
& $if$ \; \lambda(t) \neq 0, \\
a_{n-1} \lambda'(t) = -1,
& $if$ \; \lambda(t) = 0.
\endcases
$$

From Remark 4.1 we see that the factors of $\lambda'(t)$ in above equations vanish
iff the zero $\lambda(t)$ of $\phi_n^t$ is multiple. But, the right hand side can not
be null (for $\lambda(t)\neq0$ it can not be $\varphi_{n-1}^*(\lambda(t))=0$ since,
from (1.4), it would imply $\varphi_{n-1}(\lambda(t))=0$). So, we conclude that if
$\lambda(t_0)$ is multiple, then $\lambda(t)$ is not differentiable at $t=t_0$. From
Theorem 3.1, the multiplicity of $\lambda(t)$ as a zero of $\phi_n^t$ coincides with
its algebraic multiplicity as an eigenvalue of $F_n(t)$. Therefore, if $F_n(t_0)$ has
less than $n$ distinct eigenvalues, then some of them can not be expressed as
differentiable functions at $t=t_0$ and, so, $t_0 \in S$. On the contrary, all the
points $t$ where $F_n(t)$ has $n$ distinct eigenvalues are outside $S$ [41]. So, we
find that $S$ coincides exactly with the set of points where some zeros are multiple.

Finally, notice that all the quantities with index less than $n$ do not depend on $t$
since $a_j$ is constant for $j<n$. So, just using similar proofs to those given in
Proposition 2.5, Theorem 3.1 and Proposition 4.1, we see that all these results remain
true even for $t\in\T$.
$\qed$
\enddemo

Above theorem has the following immediate consequence.

\proclaim{Theorem 4.4}
Given a finite segment of orthogonal polynomials, the set of its extensions with
multiple zeros is at most denumerable. Moreover, in the positive definite case, this
set is at most finite.
\endproclaim

In a second application of Proposition 4.1, we will study the
effect of a rotation of the Schur parameters on the zeros of
orthogonal polynomials. First of all we will discuss the rotation
of only one Schur parameter. This is equivalent to analyze the
eigenvalues of
$F_n(t)=\cF_n(a_1,a_2,\dots,a_{k-1},e^{it}a_k,a_{k+1},\dots,a_n)$
as a function of the real variable $t\in[0,2\pi]$.

\proclaim{Theorem 4.5}
Let $(\varphi_n^t)_{n\geq0}$ be the orthonormal polynomials with positive leading
coefficients associated to the Schur parameters $(a_n(t))_{n\in\N}$, where
$a_n(t)=a_n\in\C\backslash\T$ for $n \neq k$ and $a_k(t)=e^{it}a_k$,
$a_k \in \C\backslash\T$.  Then, the number of distinct zeros of $\varphi_n^t$ and
their multiplicities are constant for $t\in[0,2\pi]\backslash S$, where $S$ is at most
finite. The zeros of $\varphi_n^t$ can be expressed as differentiable functions of $t$
in $[0,2\pi]\backslash S$. If $\lambda(t)$ is one of such functions for $n \geq k$,
then
$$
K_{n-1}^t(\lambda(t),\overline{\lambda(t)^{-1}}) \lambda'(t) \! = \!\!
-i\lambda(t)^{1-k} \!\!
\bigg\{ e_{k-1} a_k(t) (\varphi_{k-1}^{t*}(\lambda(t)))^2
\! + \! e_k \overline{a_k(t)} (\varphi_k^t(\lambda(t)))^2 \bigg\}.
$$
whenever $\lambda(t) \neq 0$.
\endproclaim

\demo{Proof}
First of all, notice that
$F_n(t)=\cF_n(a_1,a_2,\dots,a_{k-1},e^{it}a_k,a_{k+1},\dots,a_n)$
is not an analytic function of $t$ considered as a complex variable.
So, if we want to apply similar arguments to those given in the proof of previous
theorem, we have to change the starting point.
Let us consider the five-diagonal matrix $\cF$ given in Definition 3.1 as a function
of $a_n, \overline a_n, \rho_n, \hat\rho_n, \; n\in\N$, taking them as independent
variables. That is, denoting $a:=(a_n)_{n\in\N}$ and given the arbitrary sequences
in $\C$ $b=(b_n)_{n\in\N}$, $c=(c_n)_{n\in\N}$, $d=(d_n)_{n\in\N}$, let us define
$$
\cG(a,b,c,d) :=
\pmatrix
-a_1 & c_1 & 0 \\
-d_1 a_2  & \!\!\!\!-b_1 a_2
& \!\!\!-c_2 a_3 & c_2 c_3  \\
\;\;\,d_1 d_2 & b_1 d_2 &
\!\!\!-b_2 a_3 & b_2 c_3  & 0 \\
& \!\!\!\!\!\!0 & \!\!\!-d_3 a_4  & \!\!\!\!-b_3 a_4
& \!\!\!\!-c_4 a_5 & c_4 c_5 \\
& & \; d_3 d_4 & b_3 d_4
& \!\!\!\!-b_4 a_5 & b_4 c_5 & 0 \\
& & & \hskip-30pt\ddots & \hskip-30pt\ddots & \hskip-30pt\ddots & \hskip-15pt\ddots
& \ddots
\endpmatrix
$$
For $w\in\C^*$, let $\cG(a(w),b(w),c,d)$, where
$a_n(w)=a_n$, $b_n(w)=b_n$ for $n \neq k$ and $a_k(w)=wa_k$, $b_k(w)=w^{-1}b_k$.
Then, from $\cG_n(a(w),b(w),c,d)$ we can recover $F_n(t)$ just choosing
$b=\overline a:= (\overline a_n)_{n\in\N}$, $c=\rho:=(\rho_n)_{n\in\N}$,
$d=\hat\rho:=(\hat\rho_n)_{n\in\N}$ and $w=e^{it}$, that is,
$F_n(t)=\cG_n(a(e^{it}),\overline{a(e^{it})},\rho,\hat\rho)$.
Since $\cG_n(a(w),b(w),c,d)$ is an analytic function of $w$ in $\C^*$, analogously to
previous theorem we conclude that the number of distinct eigenvalues and their
multiplicities are constant for $w\in\C^*\backslash S'$, where $S'$ has only a finite
number of points in any compact subset of $\C$. Moreover, these eigenvalues can be
expressed as analytic functions of $w$ in any simply connected domain of
$\C^*\backslash S'$. Therefore, the zeros of $\varphi_n^t$, which are the eigenvalues
of $F_n(t)$, are constant in number and multiplicity and can be expressed as
differentiable functions of the real variable $t$ for $t\in[0,2\pi]\backslash S$,
where $S=S'\cap[0,2\pi]$ must be finite.

The rest of the theorem follows straightforward from the application of Proposition
4.1, and the simplification of the result so obtained by using Proposition 2.3 and the
relation between standard orthogonal L-polynomials and orthogonal polynomials.
$\qed$
\enddemo

Finally, we will discuss the variation of the zeros of the orthogonal polynomials under
the simultaneous rotation of all the Schur parameters. That is, we will study the
eigenvalues of $F_n(t)=\cF(e^{it}a_1,e^{it}a_2,\dots,e^{it}a_n)$ as a function of the
real variable $t\in[0,2\pi]$.

\proclaim{Theorem 4.6}
Let $(\varphi_n^t)_{n=0}$ be the orthonormal polynomials with positive leading
coefficients associated to the Schur parameters $(a_n(t))_{n\in\N}$, where
$a_n(t)=e^{it}a_n$, $a_n \in \C\backslash\T$, for all $n$.
Then, the number of distinct zeros of $\varphi_n^t$ and their multiplicities are
constant for $t \in [0,2\pi]\backslash S$, where $S$ is at most finite. Moreover, for
$t \in [0,2\pi]\backslash S$, the non null zeros of $\varphi_n^t$ are simple and can be
expressed as differentiable functions of $t$. If $\lambda(t)$ is one of such
functions, then
$$
\lambda'(t) = i\lambda(t) {1 \over K_{n-1}(\lambda(t),\overline{\lambda(t)^{-1}})}.
$$
\endproclaim

\demo{Proof}
Again, $F_n(t)=\cF(e^{it}a_1,e^{it}a_2,\dots,e^{it}a_n)$ is not an analytic function
of $t$ considered as a complex variable. So, following the notations of previous
theorem, we consider now, for each $w\in\C^*$, the infinite matrix
$\cG(a(w),b(w),c,d)$, where $a_n(w)=wa_n$, $b_n(w)=w^{-1}b_n$ for all $n$.
Since $F_n(t)=\cG_n(a(e^{it}),\overline a(e^{it}),\rho,\hat\rho)$ and
$\cG_n(a(w),b(w),c,d)$ is an analytic function of $w$ in $\C^*$, analogously to the
proof of previous theorem we conclude that the eigenvalues of $F_n(t)$ are constant in
number and multiplicity and can be expressed as differentiable functions of the real
variable $t$ for $t\in[0,2\pi]\backslash S$, where $S$ is at most finite.

Let $\lambda(t)$ be one of such functions. Just applying
Proposition 4.1 and using the relation between standard orthogonal
L-polynomials and orthogonal polynomials, together with the fact
that $\lambda(t)$ is a zero of $\varphi_n^t$, we get that, if
$\lambda(t)\neq0$,
$$
K_{n-1}^t(\lambda(t),\overline{\lambda(t)^{-1}}) \; \lambda'(t) =
-i \left(a_1(t) + \sum_{j=1}^{n-1} I_j(t)\right),
$$
where
$$
I_j(t) \kern-2pt = \kern-2pt \lambda(t)^{-j} \kern-2pt
\bigg\{ e_{j-1} \rho_j a_{j+1}(t)
\varphi_j^{t*}(\lambda(t)) \varphi_{j-1}^{t*}(\lambda(t))
+ e_j \rho_{j+1} \overline{a_j(t)}
\varphi_j^t(\lambda(t)) \varphi_{j+1}^t(\lambda(t)) \bigg\}.
$$
From (1.4) and (1.7) we find that
$$
I_j(t) =
\lambda(t)^{-j} e_j \rho_{j+1}
\varphi_{j+1}^t(\lambda(t)) \varphi_j^{t*}(\lambda(t))
-\lambda(t)^{1-j} e_{j-1} \rho_j
\varphi_j^t(\lambda(t)) \varphi_{j-1}^{t*}(\lambda(t)),
$$
which, taking again into account that $\lambda(t)$ is a zero of $\varphi_n^t$, implies
that
$$
K_{n-1}^t(\lambda(t),\overline{\lambda(t)^{-1}}) \; \lambda'(t) =
i(\rho_1\varphi_1^t(\lambda(t)) - a_1(t)) = i\lambda(t).
$$

From above result we see using Remark 4.1 that, if $\lambda(t_0) \neq 0$ is a multiple
zero, then $\lambda(t)$ can not be differentiable at $t=t_0$. Therefore, any non
vanishing zero $\lambda(t)$ must be simple for $t\in [0,2\pi]\backslash S$.
$\qed$
\enddemo

\remark{Remark 4.2} Concerning theorems 4.5 and 4.6, we have to
remark that, since $F_n(t)$ is differentiable with respect to the
real variable $t$, its eigenvalues can be expressed as
differentiable functions $\lambda(t)$ in any interval of
$[0,2\pi]$ where $F_n(t)$ is diagonable [41]. Therefore, in any
interval of $[0,2\pi]$ where the zeros of $\varphi_n^t$ are
simple, they can be expressed as differentiable functions.

If $r(t)=|\lambda(t)|$ and $\theta(t)$ is a differentiable determination for the phase
of $\lambda(t)$, then, for $\lambda(t) \neq 0$,
$$
{\lambda'(t) \over \lambda(t)} = {r'(t) \over r(t)} + i \theta'(t).
$$
Hence, Theorem 4.6 gives the following meaning for the $n-1$-th kernel $K_{n-1}(z,y)$
associated to a sequence of orthogonal polynomials on $\T$:
if $\lambda$ is a non null simple zero of the $n$-th polynomial, the real part of
$1/K_{n-1}(\lambda,\overline{\lambda^{-1})}$ measures the speed of the rotation of
$\lambda$ under rotation of the Schur parameters, while its imaginary part determines
the rapidity in the radial approach of $\lambda$ to the origin.
\endremark

\medskip

\remark{Example}
As an example, we will show the consequences of this last result in the case of the
Geronimus polynomials, defined by a constant sequence of Schur parameters $a_n=a$,
$0<|a|<1$, $n\geq1$.
In this case, it is known that the measure of orthogonality is supported in the arc
$\Delta_\alpha = \{e^{i\theta}|\alpha \leq \theta \leq 2\pi-\alpha\}$,
$\cos\alpha=1-2|a|^2$, $\alpha\in[0,\pi]$,
plus a possible mass point at $z_0=(1-a)/(1-\overline a)$ that appears iff
$\Re a > |a|^2$ [16, 18]. If we write $z_0=e^{i\psi}$, then
$\cos\psi - \cos\alpha = 2 (\Re a - |a|^2)^2 / |1-a|^2$
so, the point $z_0$ is always outside $\Delta_\alpha$, except in the case
$\Re a = |a|^2$, for which $z_0$ is an extremum $z_{\pm}=e^{\pm i\alpha}$ of the arc
$\Delta_\alpha$.

The orthonormal polynomials and their reversed are given by [16, 18]
$$
\align
& \varphi_n(z;a) = {1 \over \rho^n} (u_{n+1}-(1-a)u_n), \\
& \varphi_n^*(z;a) = {1 \over \rho^n} (u_{n+1}-(1-\overline a)zu_n),
\endalign
$$
where $\rho=\sqrt{1-|a|^2}$, $u_n = (w_1^n-w_2^n)/(w_1-w_2)$ and $w_1, w_2$
are the solutions of the quadratic equation $w^2-(z+1)w+\rho^2z=0$, that is,
$$
w_{1,2} = {1 \over 2}\left\{ {z+1 \pm \sqrt{(z+1)^2-4\rho^2z}} \right\}.
$$
Notice that $w_1+w_2=z+1$, $w_1w_2=\rho^2z$,
$(w_1-w_2)^2=(z-z_+)(z-z_-)$ and the derivatives with respect to $z$ are given by
$$
{w'_i \over w_i} = {1 \over z} {w_i-1 \over w_i-w_j}, \quad i \neq j.
$$

The zeros of $\varphi_n(z;a)$ are the solutions of $u_{n+1}=(1-a)u_n$, which implies
$$
{w_1^n \over w_2^n} = {w_2-(1-a)\over w_1-(1-a)}.
\tag 4.2
$$
If $z$ is a zero of $\varphi_n(z;a)$, (4.1) yields
$$
K_{n-1}(z,\overline{z^{-1}};a) = {z^{1-n} \over \rho^{2n}}
((1-a)-(1-\overline a)z) W(u_n,u_{n+1}).
\tag 4.8
$$
being
$$
W(u_n,u_{n+1}) = \det \pmatrix u_n & u_{n+1} \\ u'_n & u'_{n+1} \endpmatrix
$$
the Wronskian determinant of $u_n, u_{n+1}$ considered as functions of $z$.
From the definition of $u_n$ we get
$$
\align
& W(u_n,u_{n+1}) =
{1 \over (w_1-w_2)^2} \sum_{i,j=1,2} (-1)^{i+j} W(w_i^n,w_j^{n+1}) \\
& = {w_1^n w_2^n\over (w_1-w_2)^2}
\left\{ {w_1^n \over w_2^n} w'_1 + {w_2^n \over w_1^n} w'_2
- n (w_1-w_2) \left({w'_1 \over w_1}-{w'_2 \over w_2}\right) - (w'_1+w'_2) \right\}.
\endalign
$$
Using the properties of the functions $w_1$, $w_2$, together with the equation (4.2)
for the zeros, we find finally that
$$
K_{n-1}(z,\overline{z^{-1}};a) =
{ (n(z-1)+z) (1-\overline a) (z-z_0) + 2(|a|^2 - \Re a)z \over (z-z_+) (z-z_-)}.
$$

Let us consider the case $\Re a < |a|^2$, where the support of the
measure is just the arc $\Delta_\alpha$. Since $z=|z|e^{i\theta}$
is a zero of $\varphi_n(z;a)$, it has to lie on the convex hull of
the support of the measure, that is, in the subset of the open
unit disk $\D$ determined by $\cos\theta < \cos\alpha$. Thus,
$|z-1|>1-\cos\alpha=2|a|^2$,
$|z-z_0|>\cos\psi-\cos\alpha=2(|a|^2-\Re a)^2/|1-a|^2$ and
$|(z-z_+)(z-z_-)|<2(1+\cos\alpha)=4\rho^2$. Hence,
$$
|K_{n-1}(z,\overline{z^{-1}};a)| >
{(|a|^2 - \Re a) \over 2\rho^2}
\left\{ {(2n|a|^2-1) (|a|^2-\Re a) \over |1-a|} - 1 \right\}.
\tag 4.3
$$
Above inequality implies that, for $n$ big enough, $K_{n-1}(z,\overline{z^{-1}};a)$
can not vanish when $z$ is a zero of $\varphi_n(z;a)$, and, so, $\varphi_n(z;a)$ can
not have multiple zeros. More precisely, if $\Re a < |a|^2$, then
$$
n > {1 \over 2|a|^2} \left( {|1-a| \over |a|^2-\Re a} + 1 \right) \Rightarrow
\roman{the \; zeros \; of} \; \varphi_n(z;a) \; \roman{are \; simple}.
\tag 4.4
$$

We can apply now Theorem 4.6 to the perturbation $a(t)=ae^{it}$, being $a$ such that
$\Re a < |a|^2$. Since $|a|^2-\Re a = |a-1/2|^2-1/4$, we are dealing with Schur
parameters in the region of $\D$ outside the closed disk
$D_0=\{z\in\C||a-1/2|\leq1/2\}$. Without loss of generality we can suppose $a\in(0,1)$.
Given $0<t_0<t_1<2\pi$ such that $a(t_0),a(t_1)\in\D\backslash\D_0$, for
$t\in[t_0,t_1]$ it must be $a(t)\in\D\backslash\D_0$ and
$c(t):=|a(t)|^2-\Re a(t) \geq c_0:=\min\{c(t_0),c(t_1)\}$.
Therefore, from (4.4) we find that
$$
n > {1 \over 2a^2} \left( {1+a \over c_0} + 1 \right) \Rightarrow
\roman{the \; zeros \; of} \; \varphi_n(z;ae^{it}) \; \roman{are \; simple} \;
\forall t \in [t_0,t_1].
$$
Thus, from Remark 4.2 we see that, under above condition for $n$, the zeros of
$\varphi_n(z;ae^{it})$ can be expressed as differentiable functions $z(t)$ in the
interval $[t_0,t_1]$. From Theorem 4.6 and inequality (4.3) we have that
$$
\left|{z'(t) \over z(t)}\right| < C_n(a,t_0,t_1) :=
{2\rho^2 \over c_0} {1+a \over (2na^2-1) c_0 - (1+a)},
$$
and, thus, if $z(t)=|z(t)|e^{i\theta(t)}$, $\theta(t)$ differentiable in $[t_0,t_1]$,
we get the bounds
$$
||z(t_1)|-|z(t_0)||,|\theta(t_1)-\theta(t_0)| < C_n(a,t_0,t_1)|t_1-t_0|.
$$
This means that, for $n$ big enough, each zero of $\varphi_n(z;ae^{it_0})$ has a zero
of $\varphi_n(z;ae^{it_1})$ at a radial and angular distance less than
$C_n(a,t_0,t_1)|t_1-t_0|$.
Notice that
$$
C_n(a,t_0,t_1)={\rho^2 (1+a) \over a^2 c_0^2} {1 \over n}
+ O\left({1 \over n^2}\right).
$$
When $\Re a(t_0), \Re a(t_1) \leq 0$ we can choose a coefficient $C_n(a)$ independent
of $t_0,t_1$ because, then, $c_0 \geq a^2$ and, hence,
$$
C_n(a,t_0,t_1) \leq C_n(a) := {2\rho^2 \over a^2} {1+a \over 2na^4 - (1+a+a^2)},
\quad \roman{for} \;\; n > {1+a+a^2 \over 2a^4}.
$$
Notice that $C_n(a)$ is analytic for $a=1$, with $C(1)=0$. This
means that the best bounds for the location of zeros appear when
$a$ is close to 1. On the contrary, the behavior of $C_n(a)$ for
$a=0$ is singular. This suggests a more chaotic behavior for the
zeros under perturbations of the Schur parameters as far as the
last ones approach to the origin, and a slower and more regular
variation when they are close to the unit circle.
\endremark

\medskip

These are just some examples of the utility of the five-diagonal
matrix representation given for orthogonal polynomials on the unit
circle. A similar discussion to the one given here is possible for
para-orthogonal polynomials, but this together with some
applications related to the orthogonality measure will be
developed in a separated paper [9]. There it will be discussed too
an operator theoretic approach to the study of the orthogonality
measure based on the five-diagonal matrix representation obtained
for the multiplication operator.

\head
Acknowledgements
\endhead

The work of the first and second authors was supported by Direcci\'on General de
Ense\~nanza Superior (DGES) of Spain under grant PB 98-1615.
The work of the last author was supported by CAI, ``Programa Europa de Ayudas a la
Investigaci\'on".

The authors are very grateful to Professor F. Marcell\'an for his remarks and useful
suggestions, as well as to Professors E. K. Ifantis, C. G. Kokologiannaki and P. D.
Siafarikas for fruitful discussions.

\Refs

\item{[1]}  Akhiezer, N.I.
The Classical Moment Problem.
Oliver and Boyd. London, 1965.

\item{[2]}  Akhiezer, N.I.; Krein, M.G.
Some questions in the theory of moments.
Trans. Math. Mono., vol. 2, AMS. Providence, RI, 1962; Kharkov, 1938 [in russian].

\item{[3]}  Alfaro, M.
Una expresi\'on de los polinomios ortogonales sobre la circunferencia unidad.
Actas III J.M.H.L. (Sevilla, 1974), 2 (1982) 1-8.

\item{[4]}  Ambroladze, M.
On exceptional sets of asymptotics relations for general orthogonal polynomials.
J. Approx. Theory 82 (1995) 257-273.

\item{[5]}  Barrios, D.; L\'opez, G.
Ratio asymptotics for orthogonal polynomials on arcs of the unit circle.
Constr. Approx. 15 (1999) 1-31.

\item{[6]}  Barrios, D.; L\'opez, G.; Mart\'{\i}nez, A.;
Torrano, E.
On the domain of convergence and poles of $J$-fractions.
J. Approx. Theory 93 (1998) 177-200.

\item{[7]}  Barrios, D.; L\'opez, G.; Torrano, E.
Location of zeros and asymptotics of polynomials satisfying three-term recurrence
relations with complex coefficients.
Russian Acad. Sci. Sb. Math. 80 (1995) 309-333.

\item{[8]}  Bultheel, A.; Gonz\'alez-Vera, P.; Hendriksen, E.;
Nj\"astad, O.
Orthogonal Rational Functions.
University Press. Cambridge, 1999.

\item{[9]}  Cantero, M.J.; Moral, L.; Vel\'azquez, L.
Five-diagonal matrices, para-orthogonal polynomials and measures on the unit circle.
In preparation.

\item{[10]}  Chihara, T.S.
An Introduction to Orthogonal Polynomials.
Gordon and Breach. New York, 1978.

\item{[11]}  Dehesa, J.S.
The asymptotical spectrum of Jacobi matrices.
J. Comput. Appl. Math. 3 (1977) 167-171.

\item{[12]}  Delsarte, P.; Genin, Y.
Tridiagonal approach to the algebraic environment of Toeplitz matrices I. Basic
results.
SIAM J. Matrix Anal. Appl. 12 (1991) 220-238.

\item{[13]}  Delsarte, P.; Genin, Y.
Tridiagonal approach to the algebraic environment of Toeplitz matrices II. Zero
and eigenvalue problems.
SIAM J. Matrix Anal. Appl. 12 (1991) 432-448.

\item{[14]}  Dombrowski, J.
Orthogonal polynomials and functional analysis.
Orthogonal polynomials: Theory and Practice (P. Nevai, Ed.), pp. 147--161.
NATO-ASI Series C, vol. 294. Kluwer. Dordrecht, 1990.

\item{[15]}  Garc\'{\i}a, P.; Marcell\'an, F.
On zeros of regular orthogonal polynomials on the unit circle.
Ann. Polon. Math. 58 (1993) 287-298.

\item{[16]}  Geronimus, Ya.L.
Orthogonal Polynomials.
Consultants Bureau. New York, 1961.

\item{[17]}  Godoy, E.; Marcell\'an, F.
Orthogonal polynomials on the unit circle: distribution of zeros.
J. Comput. Appl. Math. 37 (1991) 195-208.

\item{[18]}  Golinskii, L.; Nevai, P.; Van Assche, W.
Perturbation of orthogonal polynomials on an arc of the unit circle.
J. Approx. Theory 83 (1995) 392-422.

\item{[19]}  Golub, G.H.; Van Loan, C.F.
Matrix  Computations,
The John Hopkins University Press, 3rd ed. Baltimore, 1996.

\item{[20]}  Hendriksen, E.
Moment methods in two point Pad\'e approximation.
J. Approx. Theory 40 (1984) 313-326.

\item{[21]}  Hendriksen, E.; Nijhuis, C.
Laurent-Jacobi matrices and the strong Hamburger moment problem.
Acta Applicandae Mathematicae 61 (2000) 119-132.

\item{[22]}  Hendriksen, E.; Van Rossum, H.
Moment methods in Pad\'e approximation.
J. Approx. Theory 35 (1982) 250-263.

\item{[23]}  Hendriksen, E.; Van Rossum, H.
Orthogonal Laurent polynomials.
Indag. Math. 48 (1986) 17-36.

\item{[24]}  Horn, R; Johnson, Ch.
Matrix Analysis.
Cambridge University Press. Cambridge, 1985.

\item{[25]}  Ifantis, E.K.
A theorem concerning differentiability of eigenvectors and eigenvalues with some
applications.
Appl. Anal. 28 (1988) 257-283.

\item{[26]}  Ifantis, E.K.
Concavity and convexity of eigenvalues.
Appl. Anal. 41 (1991) 209-220.

\item{[27]}  Ifantis, E.K.; Kokologiannaki, C.G.; Siafarikas, P. D.
Newton sum rules and monotonicity properties of the zeros of scaled co-recursive
associated polynomials.
Methods. Appl. Anal. 3 (1996) 486-497.

\item{[28]}  Ifantis, E.K.; Panagopoulos, P.N.
On the zeros of a class of polynomials defined by a three term recurrence
relation.
J. Math. Anal. Appl. 182 (1994) 361-370.

\item{[29]}  Ifantis, E.K.; Siafarikas, P.D.
An alternative proof of a theorem of Stieljes and related results.
J. Comput. Appl. Math. 65 (1995) 165-172.

\item{[30]}  Ifantis, E.K.; Siafarikas, P.D.
Perturbation of the coefficients in the recurrence relation of a class of polynomials.
J. Comput. Appl. Math. 57 (1995) 163-170.

\item{[31]}  Ismail, M.E.H.
The variation of zeros of certain orthogonal polynomials.
Adv. in Appl. Math. 8 (1987) 111-119.

\item{[32]}  Ismail, M.E.H.
Monotonicity of zeros of orthogonal polynomials.
$q$-Series and Partitions (D. Stanton, Ed.), pp. 177--190.
IMA Vol. Math. Appl., vol. 18. Springer. New York, 1989.

\item{[33]} Ismail, M.E.H.; Muldoon, M.E.
A discrete approach to the monotonicity of zeros of orthogonal polynomials.
Trans. Amer. Math. Soc. 323 (1991) 65-78.

\item{[34]}  Ismail, M.E.H.; Zhang, R.
On the Hellmann-Feynman theorem and the variation of zeros of certain special
functions.
Adv. in Appl. Math. 9 (1988) 439-446.

\item{[35]}  Ismail, M.E.H.; Zhang, R.
The Hellmann-Feynman theorem and the variation of zeros of special functions.
Ramanujan International Symposium on Analysis (N.K. Thakare et al., Eds.),
pp. 151--183.
McMillan of India. New Delhi, 1989.

\item{[36]}  Jones, W.B.; Nj\"astad, O.
Applications of Szeg\H o polynomials to digital signal processing,
Rocky Mountain J. Math. 21 (1991) 387-436.

\item{[37]}  Jones, W.B.; Nj\"astad, O.
Orthogonal Laurent polynomials and strong moment theory: a survey.
J. Comput. Appl. Math. 105 (1999) 51-91.

\item{[38]}  Jones, W.B.; Nj\"astad, O.; Thron, W.J.; Waadeland, H.
Szeg\H o polynomials applied to frequency analysis.
Comput. Appl. Math. 46 (1993) 217-228.

\item{[39]}  Jones, W. B.; Thron, W. J.
Survey of continued fraction methods of solving moment problems.
Analytic Theory of Continued Fractions.
Lecture Notes in Math., vol. 932. Springer. Berlin, 1981.

\item{[40]}  Jones, W. B.; Thron, W.J.; Nj\"astad, O.
Orthogonal Laurent polynomials and the strong Hamburger moment problem.
J. Math. Anal. Appl. 98 (1984) 528-554.

\item{[41]}  Kato, T.
Perturbation Theory for Linear Operators.
Springer. New York, 1966.

\item{[42]}  M\'at\'e, A.; Nevai, P.
Eigenvalues of finite band-with Hilbert space operators and their application to
orthogonal poluynomials.
Can. J. Math. 41 (1989) 106-122.

\item{[43]}  M\'at\'e, A.; Nevai, P.; Van Assche, W.
The support of measures associated with orthogonal polynomials and the spectra
of the related self-adjoint operators.
Rocky Mountain J. Math. 21 (1991) 501-527.

\item{[44]}  Mhaskar, H.N.;  Saff, E.B.
On the distribution of zeros of polynomials orthogonal on the unit circle.
J. Approx. Theory. 63 (1990) 30-38.

\item{[45]}  Muldoon, M.E.
Properties of zeros of orthogonal polynomials and related functions.
J. Comput. Appl. Math. 48 (1993) 167-186.

\item{[46]}  Nevai, P.
Orthogonal polynomials, recurrences, Jacobi matrices and measures.
Progress in Approximation Theory (A.A. Gonchar and E.B. Saff, Eds.), pp. 79--104.
Springer Ser. Comput. Math., vol. 19. Springer. New York, 1992.

\item{[47]}  Nevai, P.; Totik, V.
Orthogonal polynomials and their zeros.
Acta Scient. Math. (Szeged) 53 (1-2) (1989) 99-104.

\item{[48]}  Pan, K.
Asymptotics for Szeg\H o polynomials associated with Wiener-Levinson filters.
J. Comput. Appl. Math. 46 (1993) 387-394.

\item{[49]}  Pan, K.; Saff, E.B.
Asymptotics for zeros of Szeg\H o polynomials associated with trigonometric
polynomials signals.
J. Approx. Theory 71 (1992) 239-251.

\item{[50]}  Saff, E.B.
Orthogonal polynomials from a complex perspective,
Orthogonal Polynomials: Theory and Practice (P. Nevai, Ed.), pp. 363--393.
NATO-ASI Series C, vol. 294. Kluwer. Dordrecht, 1990.

\item{[51]}  Saff, E.B.; Totik, V.
What parts of a measure support attracts zeros of the corresponding orthogonal
polynomials?.
Proc. AMS 114 (1992) 185-190.

\item{[52]}  Siafarikas, P.D.
Inequalities for the zeros of the associated ultraspherical polynomials,
Math. Ineq. Appl. 2 (1999) 233-241.

\item{[53]}  Stone, M.H.
Linear Transformations in Hilbert Space.
AMS. Providence, RI, 1932.

\item{[54]}  Szeg\H o, G.
Orthogonal Polynomials.
AMS Colloq. Publ., vol. 23, AMS, 4th ed. Providence, RI, 1975.

\item{[55]}  Torrano, E.
Interpretaci\'on matricial de los polinomios ortogonales en el caso complejo,
Doctoral Dissertation, Universidad de Santander, 1987.

\endRefs

\enddocument